\documentclass[journal,draftcls,onecolumn,letterpaper,twoside,11pt]{IEEEtran}  

\usepackage[margin=1in]{geometry}

\usepackage{amssymb, mathrsfs} 
\usepackage[cmex10]{amsmath}
\usepackage{cite}
\RequirePackage[colorlinks,citecolor=blue,urlcolor=blue]{hyperref}

\usepackage{bbm}
\usepackage{epstopdf}
\usepackage{graphicx} 

\usepackage{float,latexsym,amsfonts,amstext,times}

\usepackage{color}

\newtheorem{theorem}{Theorem}[section]
\newtheorem{definition}{Definition}

\newtheorem{lemma}{Lemma}[section]
\newtheorem{remark}{Remark}[section]   

  \renewcommand{\thelemma}{\arabic{section}.\arabic{lemma}}

\newcommand{\bN}{\mathbb{N}}
\newcommand{\bR}{\mathbb{R}}

\newcommand{\um}{\underline{m}}
\newcommand{\om}{\overline{m}}

\newcommand{\cA}{\mathcal{A}}

\newcommand{\cB}{\mathcal{B}}

\newcommand{\cL}{\mathcal{L}}
\newcommand{\cF}{\mathcal{F}}

\newcommand{\cN}{\mathcal{N}}

\newcommand{\cP}{\mathcal{P}}

\newcommand{\ccab}{\mathcal{C}_{\alpha, \beta}}

\newcommand{\Fc}{{ \mathscr{F}}}

\newcommand{\Fct}{\mathscr{F}_{t}}

\newcommand{\calo}{{\mathcal{O}}}
\newcommand{\Exp}{{\sf E}}

\newcommand{\Pro}{{\sf P}}

\newcommand{\Hyp}{{\sf H}}
\newcommand{\bY}{{\mathbf{Y}}}

\newcommand{\ind}[1]{\mathbbm{1}_{\{#1\}}}  

\newcommand{\xra}{\xrightarrow} 

\newcommand{\abs}[1]{\left\vert#1\right\vert}
\newcommand{\set}[1]{\left\{#1\right\}}

\newcommand{\brc}[1]{\left(#1\right)}
\newcommand{\brcs}[1]{\left[#1\right]}

\newcommand{\Nbb}{\mathbb{N}}
\newcommand{\Xb}{\mathbf{X}}

\begin{document}

\title{Multichannel Sequential Detection---Part I: Non-i.i.d.\ Data \thanks{The work of the first author was supported by the US National Science Foundation under Grant CCF 1514245, as well as by a collaboration grant from the Simons Foundation. The work of the second author was supported in part by the U.S.\ Air Force Office of Scientific Research under MURI grant FA9550-10-1-0569, by the U.S.\  Defense Advanced Research Projects Agency under grant W911NF-12-1-0034 and by the U.S.\ Army Research Office under grant W911NF-14-1-0246 }}

\author{Georgios~Fellouris\thanks{G. Fellouris is with the Statistics Department, University of Illinois, IL 61822, USA, e-mail: fellouri@illinois.edu.}
 and Alexander~G.~ Tartakovsky, \IEEEmembership{Senior~Member,~IEEE}  
 \thanks{ A.G. Tartakovsky is at Somers, CT 06071, USA,  e-mail: alexg.tartakovsky@gmail.com .}
\thanks{Manuscript received ~~~~, 2016; revised ~~~~, .}
}

\markboth{}
{Fellouris and Tartakovsky: Multichannel Sequential Detection}

\maketitle

\begin{abstract}
We consider the problem of sequential signal detection in a multichannel system where the number and location of signals is {\em a priori} unknown. We  assume that the data in each channel are 
sequentially observed and follow a general non-i.i.d.\ stochastic model.  Under the assumption that the local  log-likelihood ratio processes in the channels  converge $r$-completely 
to positive and finite numbers,   we establish the asymptotic  optimality  of a generalized sequential likelihood ratio test and a mixture-based sequential likelihood ratio test. 
Specifically, we show that  both tests  minimize  the first $r$ moments of the stopping time distribution  asymptotically as the probabilities of  false alarm and  missed detection approach zero.  
Moreover, we show that both tests asymptotically minimize all moments of the stopping time distribution when the local log-likelihood ratio processes have independent increments and simply 
obey the Strong Law of Large Numbers. This extends a result previously known in the case of i.i.d. observations when only one channel is affected.  We illustrate the general detection theory 
using several practical examples, including the detection of signals in Gaussian hidden Markov models, white Gaussian noises with unknown intensity, and testing of the first-order 
autoregression's correlation coefficient. Finally, we illustrate the feasibility of both sequential tests when assuming an upper and a lower bound on the  number of  signals and compare their 
non-asymptotic performance using a simulation study. 
 \end{abstract}


\section{Introduction}\label{sec:Intro}

\IEEEPARstart{Q}{uick} signal detection in multichannel systems is widely applicable. For example, in the medical sphere, decision-makers must quickly detect an 
epidemic present in only a fraction of hospitals and other sources of data \cite{chang,  bock,tsu}. In environmental monitoring where a large number of sensors cover a given area, decision-makers seek to detect an anomalous behavior, such as the presence of hazardous materials or intruders, that only a fraction of sensors typically capture \cite{fie, mad}. In military defense applications, 
there is a need to detect an unknown number of targets in noisy observations obtained by radars, sonars or optical  
sensors that are typically multichannel in range, velocity and space \cite{Bakutetal-book63,Tartakovsky&Brown-IEEEAES08}. In cyber security, there is a need to rapidly detect and localize malicious 
activity, such as distributed denial-of-service attacks, typically in multiple data streams \cite{szor, Tartakovsky-Cybersecurity14, Tartakovskyetal-SM06,Tartakovskyetal-IEEESP06}. In genomic applications, 
there is a need to determine intervals of copy number variations, which are short and sparse, in multiple DNA sequences \cite{Siegmund2013}.

Motivated by these and other applications,  we consider a general  sequential detection problem  where  observations are acquired \textit{sequentially} in a number of  data streams. 
The goal is to  \textit{quickly} detect   the presence of a signal while controlling the probabilities of  false alarms (type-I error) and missed detection (type-II error) below user-specified levels.  
Two scenarios are of particular interest for applications. The first is when a single signal with an unknown location is distributed over a relatively small number of channels. For example, this may be the case 
when detecting an extended target with an unknown location in a sequence of images produced by a very high-resolution sensor. 
Following the terminology of Siegmund \cite{Siegmund2013}, we call this the ``structured'' case, 
since there is a certain geometrical structure we can know at least approximately.  A different, completely ``unstructured'' scenario is when an unknown number 
of ``point'' signals affect the channels.   For example, in many target detection applications, an unknown number of point targets appear in different channels
(or data streams), and it is unknown in which  channels the signals will appear \cite{TartakovskySQA2013}.
The multistream sequential detection problem is well-studied only in the case of a single point signal present in one 
(unknown) data stream~\cite{TLY-IEEEIT03}. However, as mentioned above, in many applications, a signal (or signals) can affect multiple data streams (e.g., when detecting an unknown number of targets in 
multichannel sensor systems). In fact, the affected subset could be completely unknown (unknown number of signals), or known partially (e.g., knowing its size or an upper bound on 
its size such as a known maximal number of signals that can appear).  

To our knowledge, this version of the sequential  multichannel detection problem has not yet been studied, although it  has  recently  received significant interest in the related sequential change detection problem; see,  e.g., \cite{felsok,Mei-B2010, Xie&Siegmund-AS13}. All these  works  focus on the case of independent and identically distributed (i.i.d.) observations in the channels. 
On the contrary, our goal  is to develop  a general  asymptotic optimality theory \textit{without assuming i.i.d. observations in the channels}. Assuming a very general non-i.i.d.\ model, we focus on two  multichannel  sequential  tests,  the Generalized Sequential Likelihood Ratio Test (\mbox{G-SLRT}) and the  Mixture Sequential Likelihood Ratio Test (\mbox{M-SLRT}), which are based on the maximum and average  likelihood ratio  over all possibly affected subsets respectively. We impose  minimal conditions on the structure of the observations in channels,  postulating only a certain asymptotic stability of the corresponding  log-likelihood ratio statistics. Specifically, we  assume that the suitably normalized log-likelihood ratios in channels almost surely converge to positive and finite numbers, which can be viewed as local limiting Kullback--Leibler information numbers.  We additionally show that if the local log-likelihood ratios also have independent increments, both  the \mbox{G-SLRT} and  the \mbox{M-SLRT}  minimize asymptotically not only the expected  sample size  but  also every  moment of the sample size distribution as the probabilities of errors vanish. Thus, we extend a result previously shown only in the case of i.i.d. observations and in the special case of a single affected stream \cite{TLY-IEEEIT03}. In the general case where the  local  log-likelihood ratios do not have independent increments,  we require a certain rate of convergence in the Strong Law of Large Numbers, which is expressed in the form of $r$-complete convergence (cf.~\cite[Ch 2]{TNB_book2014}).  Under this condition, we prove that both the \mbox{G-SLRT} and  the \mbox{M-SLRT}   asymptotically minimize the first $r$ moments of the sample size distribution.    The   $r$-complete convergence condition is a relaxation of the $r$-quick convergence condition used in \cite{TLY-IEEEIT03} (in the special case of  detecting a single signal in a multichannel system).  However, its main advantage is that it is much easier to verify in practice.  Finally, we show that  both the \mbox{G-SLRT} and  the \mbox{M-SLRT} are computationally feasible, even with a large number of channels,  when  we have an upper and a lower bound on the number of signals, a general set-up that includes cases of complete ignorance as well as cases where the size of the affected subset is known. 

The rest of the paper is organized  as follows.  In Section \ref{sec:form}, we present a mathematical formulation of the multistream sequential detection problem.  In Section \ref{sec:ALB}, 
we  obtain asymptotic lower bounds for the optimal operating characteristics. We use these lower bounds in the following sections to establish asymptotic optimality properties of the proposed 
multichannel sequential
tests.  In Section~\ref{sec:GSLRT}, we introduce the \mbox{G-SLRT} and  establish its  first-order asymptotic optimality properties  with respect to an arbitrary class of possibly affected subsets.  In Section \ref{s:WSLRT}, we introduce the \mbox{M-SLRT} and show that it has similar  asymptotic optimality properties. In Section \ref{sec:Examples}, we apply our results in the context of several stochastic models.  
In Section \ref{sec:Simulations}, we   compare the non-asymptotic performance of these two tests in a simulation study.  We conclude in  Section \ref{sec:Conclusions}. Lengthy and technical proofs are presented in the Appendix. 

Higher order asymptotic optimality properties of the test procedures,  higher order approximations for the expected sample size up to a vanishing term, and asymptotic approximations for the error probabilities  will be established  in the companion paper \cite{FTPart2}. These results are based on nonlinear renewal theory. In the companion paper, we will also present simulation results which allow us to evaluate accuracy of  the obtained approximations not only for small but also for moderate error probabilities.

\section{Problem Formulation} \label{sec:form}

Suppose  that observations are sequentially acquired over time in $K$ distinct sources, to which we refer as  channels or streams. The observations  in the $k^{\rm th}$ data stream correspond to a  realization of a discrete-time stochastic process $X^{k}:=\{X^{k}_{t}\}_{t \in \Nbb}$, where $1 \leq k \leq K$ and  $\Nbb=\{1,2, \dots\}$.   Let $\Pro^{k}$  stand for the distribution of $X^k$ and $\Pro$ for the distribution of  $\Xb=(X^{1}, \ldots, X^{K})$. Throughout the paper it is assumed that the observations from 
different channels are independent,  i.e.,   $\Pro= \Pro^{1} \times \cdots \times \Pro^{K}.$ Moreover,  for each channel there are only two possibilities,   ``noise only'' or ``signal and noise'', so that:
$$\Pro^{k} = \begin{cases}
 \Pro_{0}^{k} , \quad \text{when there is ``noise only'' in channel $k$}\\
  \Pro^{k}_{1} , \quad \text{when there is ``signal and noise'' in channel $k$}
  \end{cases}.
 $$  
Here,  $\Pro_{0}^{k}$ and $\Pro_{1}^{k}$ are distinct  probability measures on the canonical space of $X^{k}$, which are 
 mutually absolutely continuous    when restricted to the $\sigma$-algebra  $\Fc^{k}_{t}:=\sigma(X_{s}^{k}; s=1, \ldots, t)$, for any  $t \in \Nbb$.    
 Let   $\Lambda^{k}_{t}$ be the Radon--Nikod\'{y}m derivative (likelihood ratio) 
of $\Pro_{1}^{k}$ with respect to  $\Pro_{0}^{k}$ given $\Fct^{k}$ and let $Z_{t}^{k}$ be the corresponding log-likelihood ratio, i.e.,  
\begin{equation} \label{LRk}
\Lambda_{t}^{k}:= \frac{\text{d} \Pro_{1}^{k}}{\text{d} \Pro_{0}^{k}}\Big|_{\Fc^{k}_{t}} \quad \text{and} \quad Z_{t}^{k}:=\log \Lambda_{t}^{k}.
\end{equation}
 Let    $\Hyp_{0}$ be the global null  hypothesis, according to which there is  ``noise only'' in \textit{all} data streams,  and  let $\Hyp^{\cA}$ be the hypothesis according to which  a signal is present in  the subset of channels $\cA \subset \{1, \ldots, K\}$ and is absent    outside of   this subset, i.e.,  
\begin{align*}
\Hyp_{0} : & \; \Pro^k= \Pro_{0}^{k} \quad \forall \; 1 \leq k \leq K, \\
\Hyp^{\cA}: & \;
\Pro^{k}  =  \begin{cases}
\Pro_{0}^{k}  & \; \text{when} \; k \notin \cA  \\	
\Pro_{1}^{k}  & \; \text{when }  \; k \in \cA 
\end{cases} .
\end{align*}
Let $\Pro_{0}$ and $\Pro^{\cA}$ be the distributions of $\Xb$ under $\Hyp_{0}$ and $\Hyp^{\cA}$, respectively, and let $\Exp_{0}$ and $\Exp^{\cA}$ be the corresponding expectations. Let  $\Lambda_{t}^{\cA}$ be the likelihood ratio of $\Hyp^{\cA}$ against $\Hyp_{0}$ given the available information from all channels up to time  $t$,
$\Fc_{t}:=  \sigma(X_{s}; 1 \leq s \leq t),$
and  let $Z_{t}^{\cA}$  be the corresponding log-likelihood ratio (LLR). Due to the assumption of independence of channels we have  
\begin{equation} \label{LR}
\Lambda_{t}^{\cA} =\frac{\text{d} \Pro^{\cA}}{\text{d} \Pro_{0}}\Big|_{\Fc_{t}}= \prod_{k \in \cA} \Lambda_{t}^{k} \quad \text{and} \quad 
 Z_{t}^{\cA} = \log  \Lambda_{t}^{\cA}= \sum_{k \in \cA} Z_{t}^{k}.
\end{equation}
Let  $\cP$ be a  class of subsets of channels, i.e.,  a family of subsets  of $\{1, \ldots, K\}$. We want to test  the global null hypothesis $\Hyp_{0}: \Pro=\Pro_{0}$   against the alternative hypothesis 
\begin{align} \label{testing}
  \Hyp_{1} &= \bigcup_{\cA \in \cP} \Hyp^{\cA} : \;  \Pro \in \{\Pro^{\cA}\}_{\cA \in \cP},
\end{align}
according to which a signal is present in a subset of channels that belongs  to class $\cP$. Thus, class $\cP$  incorporates prior information that may be available regarding the signal.   For example,   if we have an upper $(\om \leq K)$ and a lower $(\um \geq 1)$ bound on the size of the affected subset, then we write   $\cP=\cP_{\um,\om}$, where  
\begin{align} \label{upperlowerclass}
\cP_{\um,\om} &=\{\cA \subset \{1, \ldots, K\} : \um \leq  |\cA| \leq \om\}.
\end{align}
In particular,  when we know that \textit{exactly} $m$ channels can be affected,  we write  $\cP=\cP_m$, whereas  when we know  that \textit{at most} $m$ channels can be affected, we write   $\cP=\overline{\cP}_m$, where
\begin{align*} 
\cP_m &\equiv  \cP_{m, m}= \{\cA \subset \{1, \ldots, K\} : |\cA|=m\} , \\
\overline{\cP}_m &\equiv \cP_{1,m}= \{\cA \subset \{1, \ldots, K\} : 1 \leq |\cA| \leq m\}.  
\end{align*}
Note that when we do not have any prior knowledge regarding the affected subset, which can be thought of as the most difficult case for the signal detection problem,  then $\cP=\overline{\cP}_{K}$. In what follows, we generally  assume  an arbitrary class $\cP$ of possibly affected subsets unless otherwise specified.

In this work,  we want to distinguish between  $\Hyp_{0}$ and $\Hyp_{1}$ as soon as possible. Thus, we focus on sequential tests.  We say that  the pair $(\tau,d)$  is a   \textit{sequential test}  if $\tau$ is an $\{\Fc_{t}\}$-stopping time  and $d\in \{0,1\}$ is  a binary, $\Fc_{\tau}$-measurable random variable (terminal decision) such  that  $\{d=i\}=\{\tau<\infty, \Hyp_{i} \; \text{is selected} \}$,  $i=0,1$.  We are interested in sequential tests that belong to the class
\begin{equation} \label{ccab}
\ccab(\cP):= \{ (\tau,d): \Pro_{0}(d=1) \leq \alpha ~~ \text{and} ~~ \max_{\cA \in \cP} \Pro^{\cA}(d=0) \leq \beta\},
 \end{equation}
i.e., they  control the type-I (false alarm) and type-II (missed detection) error probabilities  below  $\alpha$ and $\beta$ respectively, where $\alpha,
\beta \in (0,1)$ are arbitrary, user-specified levels.  A sequential test   $(\tau^{*},d^{*})$ in  $\ccab(\cP)$  minimizes asymptotically  as $\alpha, \beta \rightarrow 0$  the first $r$ moments of the stopping time distribution under every possible scenario if 
\begin{equation} \label{r_opt}
\Exp_{0}[(\tau^{*})^q] \sim \inf\limits_{(\tau,d) \in \ccab(\cP)} \Exp_{0}[\tau^q] \quad  \text{and} \quad
\Exp^{\cA}[(\tau^{*})^q] \sim  \inf\limits_{(\tau,d) \in \ccab(\cP)} \Exp^{\cA}[\tau^q]  \quad \forall \;  \cA \in \cP
\end{equation}
for every  integer $1  \leq q  \leq r$ and $r \in \bN$. (Hereafter, we use the standard notation $x_\alpha \sim y_\alpha$ as $\alpha\to0$, which means that $\lim_{\alpha\to0} (x_\alpha/y_\alpha)=1$). Our goal is to design  sequential tests that enjoy such  asymptotic optimality properties when we have  general stochastic models for the observations in the channels. To this end, we first obtain  asymptotic lower bounds on the moments of the stopping time distribution for sequential tests in  the class $\ccab(\cP)$, and then show that these bounds are attained asymptotically  by the proposed sequential tests.

\section{Asymptotic Lower Bounds on the Optimal Performance}\label{sec:ALB}

Since we want our  analysis to allow for very general non-i.i.d.\ models for the observations in the channels, we start by imposing a minimal condition on the structure of the observations,  which guarantees only a stability  of the local LLRs. Specifically, this stability is  guaranteed by the existence of  positive numbers  $I_{0}^{k}$ and $I_{1}^{k}$ such that  the normalized LLR process $\{ Z_t^k/t\}$   converges  almost surely (a.s.)  to $-I_0^k$ under $\Pro_0^k$ and to $I_{1}^{k}$ under $\Pro_1^k$, i.e., 
\begin{equation} \label{as}
\Pro_{1}^{k} \brc{ \frac{1}{t} Z_{t}^{k} \underset{t \rightarrow \infty} \longrightarrow  I_{1}^{k}}=1 \quad \text{and} \quad 
\Pro_{0}^{k} \brc{\frac{1}{t} Z_{t}^{k} \underset{t \rightarrow \infty} \longrightarrow  -I_{0}^{k}}=1 \quad  \forall \; 1 \leq k \leq K .
\end{equation}
In other words, we assume that each local LLR process satisfies a Strong Law of Large Numbers (SLLN). 
Obviously, this condition  implies that  
\begin{align}  \label{asA}
\Pro^{\cA} \brc{ \frac{1}{t} Z_{t}^{\cA} \underset{t \rightarrow \infty} \longrightarrow  I_{1}^{\cA}}=1 \quad \text{and} \quad 
\Pro_{0} \brc{\frac{1}{t} Z_{t}^{\cA} \underset{t \rightarrow \infty} \longrightarrow  -I_{0}^{\cA}}=1, 
\end{align}
where for any subset $\cA \subset \{1, \ldots, K\}$  we set
 \begin{equation} 
 I_{0}^{\cA} = \sum_{k \in \cA} I_{0}^{k} \quad \text{and} \quad I_{1}^{\cA}= \sum_{k \in \cA} I_{1}^{k} . \label{KL}
\end{equation}
This assumption is sufficient for establishing asymptotic lower bounds for \textit{all} moments of the stopping time distribution for sequential tests in the class $\ccab(\cP)$, which are given in the next theorem. 
 We write $\alpha_{\max} = \max(\alpha,\beta)$.
 
\begin{theorem} \label{Th:LB}
If there are  positive and finite numbers $I_0^k$ and~$I_1^k$ such that  a.s.\ convergence conditions \eqref{as} hold for every $1 \leq k \leq K$, then  for any $r \in \bN$
\begin{equation}\label{MomentsLowerbounds} 
\begin{aligned}
 \liminf_{\alpha_{\max}\to 0} \frac{\inf\limits_{(\tau,d) \in  \ccab(\cP)} \Exp_{0}[\tau^r]}{|\log\beta|^r} &  \ge \brc{\frac{1}{\min_{\cA \in \cP} I_{0}^{\cA}}}^r,
\\
\liminf_{\alpha_{\max}\to 0} \frac{\inf\limits_{(\tau,d) \in  \ccab(\cP)} \Exp^{\cA}[\tau^r]}{|\log\alpha|^r} &  \ge \brc{\frac{1}{I_{1}^{\cA}}}^r \quad \forall \;  \cA \in \cP.
\end{aligned}
\end{equation}
\end{theorem}

\begin{IEEEproof}
Fix $r\in \bN$, $\cA \in \cP$, $0<\varepsilon <1$ and   let us denote by $\ccab(\cA)$ the class of sequential tests $\ccab(\cP)$, defined in \eqref{ccab}, when $\cP=\{\cA\}$, i.e., 
\[
\ccab(\cA):= \{ (\tau,d): \Pro_{0}(d=1) \leq \alpha \; \text{and} \;  \Pro^{\cA}(d=0) \leq \beta\}.
\]
Then, for any $\alpha, \beta \in (0,1)$  we clearly have $\ccab(\cP)  \subset \ccab(\cA)$ and 
\begin{align} \label{obvious}
\inf_{(\tau,d) \in \ccab(\cP)} \Exp_{0}[\tau^r]  \geq \inf_{(\tau,d) \in \ccab(\cA)} \Exp_{0}[\tau^r] , \quad
\inf_{(\tau,d) \in \ccab(\cP)} \Exp^{\cA}[\tau^r] \geq 
\inf_{(\tau,d) \in \ccab(\cA)} \Exp^{\cA}[\tau^r].
\end{align}
Define $M_\alpha^\cA = (1-\varepsilon)  |\log \alpha|/I_1^\cA$ and $M_\beta^\cA = (1-\varepsilon)  |\log \beta|/I_0^\cA$. 
Using a quite tedious change-of-measure argument, similar to that used in the proof of  Lemma 3.4.1 in Tartakovsky~{\em et al.}~\cite{TNB_book2014}, we obtain
the inequalities 
\begin{equation}\label{ProbL2}
\begin{aligned}
\Pro^\cA\brc{\tau >  (1- \varepsilon)\frac{|\log\alpha|}{I_1^{\cA}}}  & \ge 1-\beta-\alpha^{\varepsilon^2}  - 
\Pro^\cA\brc{\max_{1\le t \le M_\alpha^\cA} Z_t^\cA \ge (1+\varepsilon) I_1^\cA M_\alpha^\cA},
\\
\Pro_0\brc{\tau >  (1- \varepsilon)\frac{|\log\beta|}{I_0^{\cA}}}  & \ge 1-\alpha-\beta^{\varepsilon^2}  - 
\Pro_0\brc{\max_{1\le t \le M_\beta^\cA} (-Z_t^\cA) \ge (1+\varepsilon) I_0^\cA M_\beta^\cA} ,
\end{aligned}
\end{equation}
that hold for an arbitrary test $(\tau,d)\in\ccab(\cA)$.  By Lemma~\ref{Lem: maxSLLN} in the Appendix, 
the a.s.\ convergence conditions \eqref{as} imply that
\begin{equation}\label{addit}
\begin{aligned} 
\Pro^\cA\brc{\max_{1\le t \le M_\alpha^\cA} Z_t^\cA \ge (1+\varepsilon) I_1^\cA M_\alpha^\cA} & \to 0 \quad \text{as}~ \alpha\to 0, \\
 \Pro_0\brc{\max_{1\le t \le M_\beta^\cA} (-Z_t^\cA) \ge (1+\varepsilon) I_0^\cA M_\beta^\cA} &\to 0 \quad \text{as}~ \beta\to 0.
  \end{aligned}
\end{equation}
From \eqref{addit} and the fact that the  right-hand sides in inequalities \eqref{ProbL2} do not depend on $(\tau,d)$ we obtain
\begin{equation}\label{ProbT}
\begin{aligned}
\lim_{\alpha_{\max}\to 0} ~ \inf_{(\tau,d) \in  \ccab(\cA)} \Pro^\cA\brc{\tau >  (1- \varepsilon)\frac{|\log\alpha|}{I_1^{\cA}}} & =1,
\\
\lim_{\alpha_{\max} \to 0} ~ \inf_{(\tau,d) \in  \ccab(\cA)} \Pro_0\brc{\tau >  (1- \varepsilon)\frac{|\log\beta|}{I_0^{\cA}}} & =1.
\end{aligned}
\end{equation}
Let us now set  $T_\alpha = \tau/\abs{\log \alpha}$. Then, from  Chebyshev's inequality we obtain
\begin{align*}
\inf_{(\tau,d)\in \ccab(\cA)} \Exp^\cA[T_\alpha^r] 
&\ge \brc{\frac{1-\varepsilon}{I_1^\cA}}^r \inf_{(\tau,d)\in \ccab(\cA)} \Pro^\cA \left(T_\alpha> (1-\varepsilon)/I_1^\cA \right) ,
\end{align*}
which together with ~\eqref{obvious} and ~\eqref{ProbT}  yields
$$
 \liminf_{\alpha_{\max}\to 0} \, \inf_{(\tau,d)\in \ccab(\cP)} \Exp^\cA[T_\alpha^r] \ge \brc{\frac{1-\varepsilon}{I_1^\cA}}^r. 
 $$
Since $\varepsilon$ is an arbitrary number in  $(0,1)$, the second inequality in \eqref{MomentsLowerbounds} follows. 

To prove the first inequality in \eqref{MomentsLowerbounds}, let $T_\beta = \tau/\abs{\log \beta}$. Then again from \eqref{obvious} and   Chebyshev's inequality  we obtain
\begin{align*}
\inf_{(\tau,d) \in \ccab(\cP)} \Exp_0[T_\beta^r] &\geq  \inf_{(\tau,d) \in \ccab(\cA)} \Exp_0[T_\beta^r] \\
&\ge \brc{\frac{1-\varepsilon}{I_0^\cA}}^r \inf_{(\tau,d)\in \ccab(\cA)} \Pro_0 \brc{T_\beta> (1-\varepsilon)/I_0^\cA} 
\end{align*}
and from \eqref{ProbT} it follows that
\[
 \liminf_{\alpha_{\max}\to 0} \; \inf_{(\tau,d)\in \ccab(\cP)}  \Exp_0[T_\beta^r] \ge  \brc{\frac{1-\varepsilon}{I_0^\cA}}^r .
 \]
But this asymptotic lower bound is valid for any $\cA\in\cP$, which implies that
\[
 \liminf_{\alpha_{\max}\to 0}  \; \inf_{(\tau,d)\in \ccab(\cP)}  \Exp_0[T_\beta^r] \ge \brc{\frac{1- \varepsilon}{\min_{\cA\in\cP}I_0^\cA}}^r .
\]
Since $\varepsilon$ is an arbitrary number in  $(0,1)$, the second inequality in \eqref{MomentsLowerbounds} follows. 
\end{IEEEproof}

\begin{remark}\label{Rem:Rem1}
A close examination of the proof of  Theorem~\ref{Th:LB} shows that it holds if  for all~$\varepsilon>0$
\begin{equation}\label{ProbMaxLLRA}
\begin{aligned}
\lim_{M\to\infty} \Pro^\cA\brc{\frac{1}{M}\max_{1\le t \le M} Z_t^\cA \ge (1+\varepsilon) I_1^\cA} & =0,
\\
\lim_{M\to\infty} \Pro_0\brc{\frac{1}{M}\max_{1\le t \le M}(- Z_t^\cA) \ge (1+\varepsilon) I_0^\cA} & = 0.
\end{aligned}
\end{equation}
As shown in Lemma~\ref{Lem: maxSLLN} in the Appendix, these conditions are guaranteed by the SLLN \eqref{asA}.
\end{remark}

\section{The Generalized Sequential Likelihood Ratio Test} \label{sec:GSLRT}

Our main goal in this section is to show that, for any given class of possibly affected subsets  $\cP$,  the asymptotic lower bounds in \eqref{MomentsLowerbounds}   are attained by the  sequential test
\begin{align} \label{gslrt}
 \hat{\tau} &= \inf\set{ t: \hat{Z}_{t}  \notin (-a,b)}, \quad 
 \hat{d}:=  \begin{cases}
1 & \; \text{when}  \quad \hat{Z}_{\hat{\tau}} \geq b   \\	
0  & \; \text{when}  \quad  \hat{Z}_{\hat{\tau}} \leq -a \\
\end{cases},
\end{align}
where $a,b>0$ are thresholds that will be selected in order to guarantee that $(\hat{\tau}, \hat{d}) \in \ccab(\cP)$  and    $\{ \hat{Z}_{t} \}$ is the maximum (generalzied) log-likelihood ratio statistic 
\begin{align} \label{glr_stat}
 \hat{Z}_{t} &= \max_{\cA \in \cP} \;   Z^{\cA}_t =   \max_{\cA \in \cP} \;  \sum_{k \in \cA} Z^{k}_t.
\end{align}
We refer to the resulting sequential test $(\hat\tau,\hat{d})$  
as the  Generalized Sequential Likelihood Ratio Test (\mbox{G-SLRT}).

\subsection{Error  Control}
 Our first task is to obtain upper bounds on the error probabilities of the \mbox{G-SLRT}, which suggest threshold values that guarantee the target error probabilities.  This is the content of the  following lemma, which does not require any assumptions on the local distributions. Let $|\cP|$ denote the cardinality of class $\cP$, i.e.,   the number of possible alternatives 
in $\cP$. Note that $|\cP|$ takes its maximum value when there is no prior information regarding the subset of affected channels ($\cP=\overline{\cP}_{K}$), in which case $|\cP|=2^{K}-1$.

\begin{lemma}\label{Lem:PEupperglrt}
For any thresholds $a, b>0$,
\begin{equation}\label{PEupper}
 \Pro_0(\hat{d} =1) \le |\cP| \, e^{-b} \quad \text{and} \quad \max_{\cA\in \cP} \Pro^\cA(\hat{d}=0) \le e^{-a}.
\end{equation}
Therefore, for any target error probabilities $\alpha, \beta \in (0,1)$, we can guarantee that  $(\hat\tau,\hat{d}) \in \ccab(\cP)$ when thresholds are selected as
\begin{equation} \label{selection}
 b=|\log (\alpha/|\cP|)| \quad \text{and} \quad a=|\log \beta|.
 \end{equation}
\end{lemma}

\begin{IEEEproof}
For any $\cA\in \cP$ we have  $Z^\cA_{\hat\tau} \leq \hat{Z}_{\hat{\tau}} \leq -a$ on $\{\hat{d}=0\}$. Therefore, by Wald's likelihood ratio identity,
\begin{align} \label{com}
\Pro^\cA(\hat{d}=0) &= \Exp_0\brcs{\exp\{Z^\cA_{\hat\tau} \} ; \hat{d} =0 }  \le e^{-a},
\end{align}
which proves the second inequality in \eqref{PEupper}.  In order to prove the first inequality  we note that  on $\{\hat{d} =1\}$ 
$$
e^{b} \leq \exp\{ \hat{Z}_{\hat{\tau}} \} =  \max_{\cA \in \cP} \Lambda_{\hat{\tau}}^{\cA}  \leq  \sum_{\cA \in \cP} \Lambda_{\hat{\tau}}^{\cA}. 
$$  
For an arbitrary  $\cB \in \cP$ we have again from Wald's likelihood ratio identity that
\[
\begin{aligned}
\Pro_0(\hat{d} =1) = \Exp^{\cB}\brcs{ \frac{1}{\Lambda^{\cB}_{\hat\tau}} \; ; \hat{d} =1} &\le  e^{-b} \,  \Exp^{\cB}\brcs{\sum_{\cA\in\cP} \frac{\Lambda^\cA_{\hat{\tau}}} {\Lambda^{\cB}_{\hat{\tau}}} \; ;  \hat{d} =1} \\
&= e^{-b} \, \sum_{\cA\in\cP} \Pro^\cA(\hat{d}=1) \le |\cP| e^{-b} .
\end{aligned}
\]
The proof is complete.
\end{IEEEproof} 

\subsection{Complete and Quick Convergence}

Asymptotic lower bounds \eqref{MomentsLowerbounds} 
were established  in Theorem \ref{Th:LB}  for any non-i.i.d.\ model that satisfies  almost sure convergence conditions \eqref{as}.  In order to show that the \mbox{G-SLRT} attains these asymptotic lower bounds, we need to  strengthen these conditions  by requiring a certain rate of  convergence. For this purpose,  it is useful to recall and clarify the notions of $r$-\textit{quick} and $r$-\textit{complete} convergence.  

\begin{definition}\label{Def:rquick}
Consider a stochastic process $(Y_t)_{t \in \bN}$ defined on a probability space $(\Omega, \cF, \Pro)$ and let $\Exp$  be the expectation that corresponds to $\Pro$. Let also $r>0$ be some positive number. 

\noindent (i) We say that $(Y_t)_{t \in \bN}$ converges {\em $r$-quickly} under $\Pro$ to a constant $I$ as $t\to\infty$ and write
\[
Y_t \xra[t \to \infty]{\Pro-r-{\text{quickly}}} I ,
\]
 if $\Exp[L(\varepsilon)]^r <\infty$ for all $\varepsilon >0$, where $L(\varepsilon) = \sup\set{t\ge 1: |Y_t-I|>\varepsilon}$ is the last time $t$ that $Y_t$ leaves the interval 
$[I-\varepsilon, I+\varepsilon]$ ($\sup\{\varnothing\}=0$). 

\noindent (ii) We say that  $(Y_t)_{t \in \bN}$ converges {\em $r$-completely} under $\Pro$ to a constant $I$ as $t\to\infty$ and write
\[
Y_t \xra[t\to\infty]{\Pro-r-{\text{completely}}} I ,
\]
 if 
 $$
 \sum_{t=1}^\infty t^{r-1} \Pro\brc{|Y_t -I| > \varepsilon} < \infty \quad \text{for all}~ \varepsilon >0.
 $$
 \end{definition}

\begin{remark}
For $r=1$, $r$-complete convergence is equivalent to  complete convergence introduced by Hsu and Robbins~\cite{Hsu&Robbins-47}. 
\end{remark}

\begin{remark}
Almost sure convergence $\Pro\left(Y_t\rightarrow q\right)=1$ is equivalent to $\Pro\brc{L(\varepsilon) < \infty} =1$ for all $\varepsilon >0$. Thus, it is implied by  $r$-quick convergence for any $r>0$ and by 
 $r$-complete convergence   for any  $r\geq1$ (due to the Borel--Cantelli lemma).
 \end{remark}

\begin{remark}
It follows from Theorem~2.4.4 in \cite{TNB_book2014} that $r$-quick convergence and $r$-complete convergence are equivalent when $\{Y_{t}\}$ is an average of i.i.d.\ random variables that have a 
 finite absolute moment of order $r+1$. In this case, these types of convergence determine a rate   of convergence in the SLLN, a topic considered in detail by Baum and Katz~\cite{Baum&Katz65}. In general, $r$-quick convergence is  somewhat stronger than $r$-complete convergence (cf.\ Lemma~2.4.1 in \cite{TNB_book2014}).  More importantly,  $r$-quick convergence is  usually  more difficult to verify in particular examples. For this reason, in the present paper, we establish asymptotic optimality of the \mbox{G-SLRT} under $r$-complete, instead of r-quick, convergence conditions. 
 \end{remark}

\subsection{Asymptotic Optimality Under $r$-complete Convergence Conditions} \label{ssec:AOGSLRT}

Our next goal is to  show that  if thresholds $a, b$ are selected according to 
\eqref{Th:LB}, then the \mbox{G-SLRT} attains the asymptotic lower bounds for moments of the sample size given  in \eqref{MomentsLowerbounds}  for every integer  $1 \leq q \leq r$ 
when the local LLRs obey  a strengthened ($r$-complete) version of the SLLN, 
\begin{equation} \label{rcompleteLLRs}
\frac{1}{t} Z_t^k \xra[t\to\infty]{\Pro_0^k-r-\text{completely}}  -I_0^k
 \quad \text{and} \quad 
 \frac{1}{t} Z_t^k \xra[t\to\infty]{\Pro_1^k-r-\text{completely}}  I_1^k ,
 \quad 1 \leq k \leq K,
\end{equation} 
i.e., assuming that  for all $\varepsilon >0$,
\begin{equation} \label{rcompleteLLRsdetailed}
\begin{aligned}
\sum_{t=1}^\infty t^{r-1} \Pro_0^k\brc{\abs{\frac{1}{t} Z_{t}^{k} +I_0^k} > \varepsilon} < \infty, \quad   \sum_{t=1}^\infty t^{r-1} \Pro_1^k\brc{\abs{\frac{1}{t} Z_{t}^{k} -I_1^k} > \varepsilon} < \infty , 
       \quad 1 \leq k \leq K.
 \end{aligned}
 \end{equation}

Before we establish the main results of this section (Theorem \ref{Th:GLRTAO}), 
 we state some auxiliary results that are necessary for the proof 
but also are of independent interest.  We start with  Lemma~\ref{rcompleteglobal}  which states that $r$-complete convergence of the local LLRs guarantees $r$-complete convergence of the cumulative LLR $Z^{\cA}$. 
The proof is given in the Appendix.

\begin{lemma} \label{rcompleteglobal}
Let $r\in \bN$. If the local $r$-complete convergence conditions  \eqref{rcompleteLLRs}  hold, then for every  $\cA \in \cP$
\begin{equation} \label{rquickLLR2_new}
\frac{1}{t} Z_t^\cA \xra[t\to\infty]{\Pro^{\cA} -r-\text{completely}}  I_1^\cA \quad \text{and} \quad 
\frac{1}{t} Z_t^{\cA} \xra[t\to\infty]{\Pro_0-r-\text{completely}}  -I_0^\cA,
\end{equation} 
where  $ I_1^\cA$ and  $I_0^\cA$ are defined in \eqref{KL}.
Moreover,
\begin{equation*}
\max_{\cA \in \cP} \, \Bigl| \frac{Z_t^{\cA}}{ I_0^\cA \, t}   +1 \Bigr| \xra[t\to\infty]{\Pro_0-r-\text{completely}}   0.
\end{equation*} 
\end{lemma}

The following theorem  provides a first-order asymptotic approximation for the moments of the \mbox{G-SLRT} stopping time for large threshold values.  These asymptotic approximations  may be useful, 
apart from  proving asymptotic optimality in Theorem~\ref{Th:GLRTAO},  for  problems with different types of constraints, for example in Bayesian settings.

\begin{theorem} \label{AEmoments}
Let $r\in \bN$. If conditions \eqref{rcompleteLLRs}  are satisfied, then  the following asymptotic approximations hold
\begin{equation}\label{AEqmomentsGLRT}
 \lim_{a_{\min} \to \infty} \frac{\Exp_{0}[\hat\tau^q]}{a^q}  =\brc{\frac{1}{\min\limits_{\cA \in \cP}   I_{0}^{\cA} } }^q, \quad  
 \lim_{a_{\min} \to \infty} \frac{\Exp^{\cA}[\hat\tau^q]}{b^q}= \brc{\frac{1}{I_{1}^{\cA} }}^q ,
\end{equation}
for every integer $1 \le q \le r$ and  $\cA\in\cP$,  where $a_{\min}=\min(a,b)$. 
\end{theorem}

\begin{IEEEproof} 
Fix $\varepsilon \in (0,1)$, $\cA \in \cP$ and set $M_b^\cA = (1-\varepsilon) b/I_1^\cA$ and $M_a^\cA = (1-\varepsilon) a/I_0^\cA$. 
Similarly to \eqref{ProbL2} we obtain
\begin{equation}\label{ProbGLRT}
\begin{aligned}
\Pro^\cA\set{\hat\tau >  (1- \varepsilon)\frac{b}{I_1^{\cA}}}  & \ge 1- \Pro^\cA(\hat{d}=0) -[\Pro_0(\hat{d} =1)]^{\varepsilon^2}  - \Pro^\cA\brc{\max_{1\le t \le M_b^\cA} Z_t^\cA \ge (1+\varepsilon) I_1^\cA M_b^\cA},
\\
\Pro_0\set{\hat\tau >  (1- \varepsilon)\frac{a}{I_0^{\cA}}}  & \ge 1-\Pro_0(\hat{d} =1) -[\Pro^\cA(\hat{d}=0)]^{\varepsilon^2}  - \Pro_0\brc{\max_{1\le t \le M_a^\cA} (-Z_t^\cA) \ge (1+\varepsilon) I_0^\cA M_a^\cA}.
\end{aligned}
\end{equation}
Combining \eqref{ProbGLRT} with \eqref{PEupper} yields
\begin{equation}\label{ProbGLRT1}
\begin{aligned}
\Pro^\cA\left( \hat\tau >  (1- \varepsilon)\frac{b}{I_1^{\cA}} \right)  & \ge 1- e^{-a} -[|\cP| e^{-b}]^{\varepsilon^2}  - \Pro^\cA\brc{\max_{1 \le t \le M_b^\cA} Z_t^\cA \ge (1+\varepsilon) I_1^\cA M_b^\cA},
\\
\Pro_0\left( \hat\tau >  (1- \varepsilon)\frac{a}{I_0^{\cA}} \right)  & \ge 1-|\cP| e^{-b} -e^{-a \varepsilon^2}  - \Pro_0\brc{\max_{1\le t \le M_a^\cA} (-Z_t^\cA) \ge (1+\varepsilon) I_0^\cA M_a^\cA}.
\end{aligned}
\end{equation}
By $r$-complete convergence conditions \eqref{rcompleteLLRs}, Lemma  \ref{rcompleteglobal}, and Lemma~\ref{Lem: maxSLLN},
\begin{equation}\label{ProbmaxZA}
\Pro^\cA\left(\max_{1\le t \le M} Z_t^\cA \ge (1+\varepsilon) I_1^\cA M \right)  \xra[M\to\infty]{} 0, \quad \Pro_0\left( \max_{1\le t \le M} (-Z_t^\cA) \ge (1+\varepsilon) I_0^\cA M \right)  \xra[M\to\infty]{} 0,
\end{equation} 
so that inequalities \eqref{ProbGLRT1} imply
\[
\Pro^\cA\left( \hat\tau >  (1- \varepsilon)\frac{b}{I_1^{\cA}} \right) \xra[a_{\min}\to\infty]{} 1 , \quad \Pro_0\left( \hat\tau >  (1- \varepsilon)\frac{a}{I_0^{\cA}} \right) \xra[a_{\min}\to\infty]{}  1.
\]
Hence,  for any $q \geq 1$,   Chebyshev's inequality yields the following asymptotic lower bounds for the moments of the  stopping time of the \mbox{G-SLRT}:
\begin{equation}\label{LBmomentsGLRT}
\Exp^{\cA}[\hat\tau^q] \ge \brc{\frac{b}{I_{1}^{\cA} }}^q(1+o(1)), \quad  \Exp_{0}[\hat\tau^q] \ge \brc{\frac{a}{\min\limits_{\cA \in \cP}   I_{0}^{\cA} } }^q (1+o(1)) \quad \text{as}~ a_{\min}\to \infty.
\end{equation}
In order to obtain  asymptotic equalities  \eqref{AEqmomentsGLRT}, it suffices to establish the asymptotic upper bounds
\begin{equation}\label{UB_GLRT}
\Exp^{\cA}[\hat\tau^r] \le \brc{\frac{b}{I_{1}^{\cA} }}^r(1+o(1)), \quad  \Exp_{0}[\hat\tau^r] \le \brc{\frac{a}{\min\limits_{\cA \in \cP}   I_{0}^{\cA} } }^r (1+o(1)) \quad \text{as}~ a_{\min}\to \infty,
\end{equation}
as \eqref{AEqmomentsGLRT} would then hold for every $1  \leq q \leq r$  with an application of H\"older's inequality.  Note also that since
\begin{equation}\label{Tandnu}
\hat\tau \le T_b^\cA := \inf\set{t: Z_t^\cA \ge b} \quad \text{and} \quad \hat\tau \le \nu_a := \inf\set{t: \min_{\cA\in \cP} (-Z_t^\cA) \ge a},
\end{equation}
 it suffices to show that 
\begin{equation}\label{UB_GLRT2}
\Exp^{\cA} \left[ (T_b^\cA )^r \right] \le \brc{\frac{b}{I_{1}^{\cA} }}^r(1+o(1)) \quad \text{as}~ b \to \infty 
\end{equation}
and
\begin{equation}\label{UB_GLRT3}
\Exp_{0}[\nu_a^r] \le \brc{\frac{a}{\min\limits_{\cA \in \cP}   I_{0}^{\cA} } }^r (1+o(1)) \quad \text{as}~ a \to \infty .
\end{equation}

Let $N_b=\lfloor b/(I_1^\cA-\varepsilon) \rfloor$ be an integer number $ \le b/(I_1^\cA-\varepsilon)$.   
We have the following chain of equalities and inequalities
\begin{align}\label{IneqExpTb}
\Exp^{\cA} \left[ (T_b^\cA )^r \right]   & = \int_0^\infty r t^{r-1} \Pro^\cA\brc{T_b^\cA > t} \, dt   \nonumber
\\
&   = r \int_0^{N_b+1} t^{r-1} \Pro^\cA\brc{T_b^\cA > t} \, dt  + r \int_{N_b+1}^{\infty} t^{r-1} \Pro^\cA\brc{T_b^\cA > t} \, dt  \nonumber
\\
& \leq  (1+N_b)^r + \sum_{\ell=1}^{\infty} \int_{N_{b}+\ell}^{N_{b}+\ell+1} r t^{r-1}  \Pro^\cA (T_b^\cA > t) \, dt \nonumber
\\
& \leq (1+N_b)^r + \sum_{\ell=1}^{\infty} \int_{N_{b}+\ell}^{N_{b}+\ell+1} r t^{r-1}  \, \Pro^\cA (T_b^\cA > N_b +\ell) \, dt \nonumber
\\
& = (1+N_b)^r + \sum_{\ell=1}^{\infty} [(N_{b}+\ell+1)^r- (N_{b}+\ell)^r ]\; \Pro^\cA (T_b^\cA > N_b +\ell) \nonumber
\\
& = (1+N_b)^r + \sum_{\ell=N_b+1}^{\infty} [(\ell+1)^r-\ell^r]  \; \Pro^\cA (T_b^\cA > \ell) \nonumber
\\
 & \le   (1+N_{b})^{r} +\sum_{\ell=N_b+1}^{\infty}   r (\ell+1)^{r-1}  \;  \Pro^\cA (T_b^\cA > \ell) \nonumber
 \\
& \le  (1+N_{b})^{r} + r 2^{r-1}  \, \sum_{\ell=N_b+1}^{\infty}  \ell^{r-1}   \Pro^\cA (T_b^\cA > \ell) . 
\end{align}
Setting  $Y_t^\cA: = t^{-1} Z_t^\cA - I_1^\cA$, we observe that for any $t \in \bN$ we have 
\begin{align*}
\Pro^\cA\brc{T_b^\cA > t}  &= \Pro^\cA\brc{ \max_{1 \leq s \leq t} Z_s^\cA < b}  \le \Pro^\cA \brc{Z_t^\cA < b}  \le   \Pro^\cA\brc{Y_t^\cA < -I_1^\cA +b/t}  .
\end{align*}
Consequently,  for any $t > N_b$ and $0<\varepsilon < I_1^\cA$ we have
\begin{equation}\label{Probineq1}
\Pro^\cA\brc{T_b^\cA > t}   \le \Pro^\cA\brc{Y_t^\cA < - \varepsilon}   \le \Pro^\cA (|Y_t^\cA| > \varepsilon).
\end{equation}
Using \eqref{Probineq1}, we conclude that
\[
 \sum_{\ell=N_b+1}^{\infty} \ell^{r-1} \, \Pro^\cA  \brc{T_b^\cA > \ell} \le \sum_{\ell=1}^{\infty} \ell^{r-1} \Pro^\cA (|Y_\ell^\cA| > \varepsilon) \equiv U_r^\cA(\varepsilon),
\]
so that 
$$
\Exp^{\cA}\brcs{T_b^\cA} \le (1+N_b)^r + r 2^{r-1}U_r^\cA(\varepsilon). 
$$
By Lemma~\ref{rcompleteglobal}, $U_r^\cA(\varepsilon)<\infty$ for all $\varepsilon >0$. Consequently, for any $0<\varepsilon < I_1^\cA$,
\begin{equation}\label{UBTA}
\Exp^{\cA} \left[ (T_b^\cA )^r \right]  \le  N_b^r (1+o(1))  =  \brc{\frac{b}{I_{1}^{\cA}-\varepsilon}}^r(1+o(1)) \quad \text{as}~ b \to \infty .
\end{equation}
Letting $\varepsilon\to 0$, we obtain asymptotic upper bound \eqref{UB_GLRT2}, which 
along with lower bound \eqref{LBmomentsGLRT} implies the second asymptotic approximation in \eqref{AEqmomentsGLRT}. 

Next, define the Markov time
\[
\tilde \nu_a = \inf\set{t \in \bN: \min_{\cA\in \cP} (-\widetilde{Z}_t^\cA) \ge \tilde{a}},  \; \text{where} \; \widetilde{Z}_t^\cA := Z_t^\cA/I_0^\cA, \quad \tilde{a} := a/\min_{\cA\in\cP} I_0^\cA.
\]
Clearly, $\nu_a \le \tilde \nu_a$, so in order to obtain upper bound \eqref{UB_GLRT3} it suffices to prove that  
this bound holds for $\Exp_{0}[\tilde \nu_a^r]$. Let $\widetilde{Y}_t^\cA= \min_{\cA\in \cP} (-\widetilde{Z}_t^\cA) / t +1$. We have
\[
\Pro_0(\tilde \nu_a > t) \le \Pro_0 \brc{\min_{\cA\in \cP} (-\widetilde{Z}_t^\cA)< \tilde{a}} = \Pro_0\brc{\widetilde{Y}_t^\cA < -1 + \tilde{a}/t} .
\]
Set  $N_a =\lfloor\tilde{a}/(1-\varepsilon)\rfloor=\lfloor a/(\min_{\cA\in\cP}I_1^\cA(1-\varepsilon))\rfloor$. Then,  for any $0<\varepsilon <1$ and  $t> N_a$, we have 
$$
\Pro_0\brc{\widetilde{Y}_t^\cA < -1 + \tilde{a}/t}  \leq \Pro_0(\widetilde{Y}_t^\cA < -\varepsilon) \leq  \Pro_0( |\widetilde{Y}_t^\cA| > \varepsilon)
$$
 and, consequently, 
\begin{equation}\label{Probineq2}
\Pro_0\brc{\tilde \nu_a > t} \le \Pro_0( |\widetilde{Y}_t^\cA| > \varepsilon). 
\end{equation}
Now, applying the same argument as above that has led to \eqref{IneqExpTb}, we obtain
\[
\Exp_0 \left[\tilde \nu_a^r \right]  \le   (1+N_{a})^{r} +r 2^{r-1}  \sum_{\ell=N_a+1}^{\infty} \ell^{r-1} \Pro_0 (\tilde \nu_a > \ell) ,
\]
which along with inequality \eqref{Probineq2} yields
\[
\Exp_0 \left[\tilde \nu_a^r\right]  \le   
 (1+N_{a})^{r} +  r 2^{r-1}\sum_{\ell=1}^{\infty}   \ell^{r-1} \Pro_0 (|\widetilde{Y}_\ell | > \varepsilon),
\]
where the last sum is finite by the $r$-complete convergence \eqref{rquickLLR2_new} (see  Lemma~\ref{AEmoments}). Therefore, for any $0<\varepsilon<1$,
\begin{equation}\label{UBnu}
\Exp_0 [\nu_a^r] \le \Exp_0 [\tilde \nu_a^r] \le N_a^r(1+o(1)) = \brc{\frac{a} {  (1-\varepsilon) \min\limits_{\cA \in \cP}   I_{0}^{\cA}  }}^r(1+o(1)) \quad \text{as}~ a \to \infty .
\end{equation}
Since $\varepsilon\in(0,1)$ is arbitrary, this implies upper bound \eqref{UB_GLRT3} and hence the second asymptotic upper bound in \eqref{UB_GLRT}.  The proof of asymptotic equalities \eqref{AEqmomentsGLRT} is complete.
\end{IEEEproof}

We are now prepared to prove the following theorem, which establishes first-order asymptotic optimality of the \mbox{G-SLRT} with respect to positive moments of the 
stopping time distribution.  

\begin{theorem}\label{Th:GLRTAO}
Consider an arbitrary class of alternatives, $\cP$, and suppose that  the thresholds $a,b$ of the \mbox{G-SLRT}  are  chosen  according to \eqref{selection}. If  the $r$-complete convergence conditions \eqref{rcompleteLLRs}  hold  for some $r\in \bN$, then for  all $1 \le q\le r$ we have   as $\alpha_{\max} \to 0$,
 \begin{equation}\label{gsprtAO0}
 \Exp_{0}[\hat\tau^q] \sim \brc{\frac{|\log \beta|}{\min\limits_{\cA \in \cP}   I_{0}^{\cA} } }^q \sim \inf_{(\tau,d) \in \ccab(\cP)} \Exp_{0}[\tau^q] 
\end{equation}
and for every $\cA \in \cP$
 \begin{equation} \label{gsprtAO1} 
 \Exp^{\cA}[\hat\tau^q] \sim \brc{\frac{|\log \alpha|}{I_{1}^{\cA} }}^q \sim \inf_{(\tau,d) \in \ccab(\cP)} \Exp^{\cA}[\tau^q] .
\end{equation}
\end{theorem}

\begin{IEEEproof}
From \eqref{PEupper} it follows that if we set $b=|\log (\alpha/|\cP|)|$ and $a=|\log \beta|$, then   $(\hat\tau,\hat{d}) \in  \ccab(\cP)$. Substituting these threshold values into asymptotic approximations \eqref{AEqmomentsGLRT}, we obtain 
 \[
 \Exp_{0}[\hat\tau]^q \sim \brc{\frac{|\log\beta|}{\min\limits_{\cA \in \cP}   I_{0}^{\cA} } }^q, \quad  \Exp^{\cA}[\hat\tau]^q \sim \brc{\frac{|\log\alpha|}{I_{1}^{\cA} }}^q \quad \text{as}~ \alpha_{\max}\to 0.
\]
Comparing with lower bounds \eqref{MomentsLowerbounds} in Theorem~\ref{Th:LB} proves \eqref{gsprtAO0}--\eqref{gsprtAO1}. 
\end{IEEEproof}

\begin{remark}
The theorem remains valid for any selection of thresholds such that  $(\hat\tau,\hat{d}) \in  \ccab(\cP)$ and $b\sim|\log \alpha|$, $a\sim |\log \beta|$ as  $\alpha_{\max} \to 0$.
\end{remark}

\begin{remark}\label{Rem:Completeinstead}
A closer examination of the proofs of  Lemma~\ref{AEmoments} and Theorem~\ref{Th:GLRTAO} 
shows that their assertions hold if the $r$-complete convergence conditions are replaced by the left-tail conditions
\[
\sum_{t=1}^\infty t^{r-1} \Pro_0^{k} \brc{-\frac{1}{t} Z_t^k < I_0^k - \varepsilon} < \infty, \quad \sum_{t=1}^\infty t^{r-1} \Pro_1^k\brc{\frac{1}{t}Z_t^k < I_1^k - \varepsilon}  
\quad \text{for all}~\varepsilon >0, ~1 \leq k \leq K 
\]
along with the SLLN in \eqref{as},  i.e.,  $\Pro_0(t^{-1}  Z_t \to -I_0^k)$, $ \Pro_1^k(t^{-1}Z_t \to  I_1^k)$ as $t\to\infty$.
 In fact, it can be shown that these conditions guarantee the uniform integrability of the sequences $\{T^{\cA}_b\}_{b>1}$ and $\{\nu_a\}_{a>1}$, 
defined in \eqref{Tandnu}, and this  can be used for an alternative proof of the theorem.
  \end{remark}

\begin{remark}
Theorem  \ref{Th:GLRTAO} was established in  \cite{TLY-IEEEIT03} in the special case where the signal  can be present in only one channel,  i.e., 
when    $\cP= \cP_{1}$, under the stronger (and harder to check) $r$-quick  convergence conditions
\begin{equation*} 
\frac{1}{t} Z_t^k \xra[t\to\infty]{\Pro_0^k-r-\text{quickly}}  -I_0^k \quad \text{and} \quad 
\frac{1}{t} Z_t^k \xra[t\to\infty]{\Pro_1^k-r-\text{quickly}}  I_1^k,  \quad 1 \leq k \leq K.
\end{equation*}
\end{remark}

\subsection{Asymptotic Optimality when the LLRs Have Independent Increments}\label{ssec:AOindep}

Let  $\ell_t^{k}=Z_{t}^k-Z_{t-1}^k$, $t \in \bN$ be the sequence of LLR increments in the $k^{\rm th}$ channel.  We now show that if each  $(\ell_t^{k})_{t \in \bN}$ is a sequence of independent, \textit{but not necessarily identically distributed}, random variables,   the asymptotic optimality properties \eqref{gsprtAO0}--\eqref{gsprtAO1} hold true for any positive integer $q$, as long as only the  a.s. conditions \eqref{as} are satisfied. To this end, we need  the following  renewal theorem, whose  proof is presented in the Appendix.

\begin{lemma} \label{independent}
Let $\xi^{k}:=(\xi_t^k)_{t \in \bN}$, $1 \leq k \leq K$ be (possibly dependent) 
sequences of random variables  on some probability space $(\Omega, \cF, \cP)$ and let $\Exp$ the corresponding expectation.  Define
the stopping time 
$$\nu(b) := \inf\left\{ t \in \bN : \min_{1\leq k \leq K} S_t^k >b \right\}; \quad  S_t^k:=\sum_{u=1}^{t} \xi_u^k.
$$ 
Suppose that  for every  $1 \leq k \leq K$  there is a positive constant $\mu_k$ such that 
$
S_t^k / t  \overset{a.s.} \longrightarrow \mu_k.
$
Then, as $b \rightarrow \infty$ we have 
$$ \frac{\nu(b)}{b} \overset{a.s.} \longrightarrow  \left(\min_{1 \leq k \leq K} \mu_k \right)^{-1}.
$$
Moreover, the convergence holds in $\cL^{r}$ for every $r>0$, if each $\xi^k$ is a sequence of  independent random variables and there is a  $\lambda \in (0,1)$ such  that 
\begin{equation} \label{condition}
\sup_{t \in \bN} \Exp \left[ \exp\{\lambda (\xi_t^k)^{-} \} \right] <\infty.
\end{equation}

\end{lemma}

The following theorem establishes a stronger asymptotic optimality property for the \mbox{G-SLRT} in the case of LLRs with independent increments.

\begin{theorem}\label{Th:GLRTAOindep_gslrt}
Let $\cP$ be an arbitrary class of  possibly affected subsets of channels
and suppose that the thresholds in the \mbox{G-SLRT}  are selected  according to \eqref{selection}.  If the LLR increments, $\{\ell_{t}^{k}\}_{t\in\bN}$,   are independent over time under 
$\Pro_0^k$ and $\Pro_1^k$ for every $1 \leq k \leq K$, then the asymptotic optimality properties \eqref{gsprtAO0}--\eqref{gsprtAO1} hold true for any $q \in \bN$, as long as the almost sure convergence conditions  \eqref{as} hold. 
\end{theorem}

\begin{IEEEproof}
By Theorem~\ref{Th:LB}, asymptotic lower bounds \eqref{MomentsLowerbounds} hold, so it suffices to show that when the thresholds in the \mbox{G-SLRT}  are selected  according to \eqref{selection},   for all $r \in \bN$ we have
\[
\limsup_{\alpha_{\max} \to 0} \frac{\Exp^\cA [\hat\tau^r]}{|\log\alpha|^r} \le \brc{\frac{1}{I_{1}^{\cA}}}^r , 
\quad \limsup_{\alpha_{\max}\to 0} \frac{\Exp_0 [\hat\tau^r]}{|\log\beta|^r} \le \brc{\frac{1}{\min_{\cA\in\cP}I_0^{\cA}}}^r. 
\]
Recall now the inequalities \eqref{Tandnu}, according to which 
\begin{equation}\label{Tandnu2}
\hat\tau \le T_b^\cA := \inf\set{t: Z_t^\cA \ge b} \quad \text{and} \quad \hat\tau \le \nu_a := \inf\set{t: \min_{\cA\in \cP} (-Z_t^\cA) \ge a}.
\end{equation}
Then  it is clear that it suffices to show that 
\[
\lim_{b \to\infty} \frac{\Exp^\cA \left[ (T_b^\cA)^r \right]}{ b^r} = \brc{\frac{1}{I_{1}^{\cA}}}^r , 
\quad \lim_{a\to\infty} \frac{\Exp_0 \left[\nu_a^r \right]}{a^r} = \brc{\frac{1}{\min_{\cA\in\cP}I_0^{\cA}}}^r. 
\]
This follows directly  from  Lemma~\ref{independent} as soon as we show that there is a  $\lambda \in (0,1)$ such that
\begin{equation} \label{cond}
\sup_{t \in \bN} \Exp^{\cA} \left[ \exp\{\lambda (\ell_t^{\cA})^{-} \}\right] <\infty,
\quad \sup_{t \in \bN} \Exp_{0} \left[ \exp\{ \lambda (-\ell_t^{\cA})^{-} \} \right] <\infty,
\end{equation}
where  $\ell^{\cA}$ are the  increments of  $Z^{\cA}$, i.e.,
$$
\ell_t^{\cA}=Z_{t}^\cA-Z_{t-1}^\cA= \sum_{k \in \cA} \ell_t^k,  \quad t \in \bN.$$ 
Indeed,  for any given $\lambda \in (0,1)$, from Jensen's inequality we have 
$$
\Exp^{\cA}[ \exp\{\lambda (\ell_t^{\cA})^{-} \}] \leq  \Exp^{\cA}[ \exp\{-\lambda \ell_t^{\cA} \} ] + 1  <  \Exp^{\cA}[ \exp\{ -\ell_t^{\cA} \}  ] ^{\lambda} +1  =2,
$$  
where the equality holds because each   $\exp\{ -\ell_t^{\cA} \} $ has mean 1  under $\Pro^{\cA}$, as a likelihood ratio. The second condition in \eqref{cond} can be verified in a similar way.
\end{IEEEproof}

\begin{remark}
The LLR increments  $(\ell_t^{\cA})_{t \in \bN}$ can be independent over time not only when the acquired observations $\{X_t\}$ are independent over time, but also  for certain models of dependent
observations that produce a sequence of LLRs with independent increments.  See, e.g., an example in Subsection~\ref{ssec:Kalman}.
\end{remark}

\begin{remark}
This result was  obtained   in \cite{TLY-IEEEIT03}  in the special case that the 
LLR increments  $(\ell_t^{k})_{t \in \bN}$ in each stream are  independent \textit{and identically distributed}  and a signal can be  present in at most one stream, i.e.,  $\cP=\cP_1$. 
\end{remark}

\subsection{Feasibility}\label{ssec:GSLRT}

The implementation of  the \mbox{G-SLRT} requires computing at each time $t$ the generalized log-likelihood ratio statistic  \eqref{glr_stat}, 
\begin{align*}
 \hat{Z}_{t} &= \max_{\cA \in \cP} \;   Z^{\cA}_t =   \max_{\cA \in \cP} \;  \sum_{k \in \cA} Z^{k}_t.
\end{align*}
 A direct computation of each $Z_t^{\cA}$ for every $\cA  \in \cP$ can be a very computationally expensive task when the cardinality of class $\cP$, $|\cP|$, is very large. However,  the computation of   $\hat{Z}_{t}$  is very easy  for a  class $\cP$  of the form  $\cP_{\um, \om}$,  which  contains all subsets of size at least $\um$ and at most $\om$. In order to see this, let  us use the following notation for the order statistics: $Z_t^{(1)} \geq \ldots \geq Z_t^{(K)}$,   i.e., $Z_t^{(1)}$ is the top local LLR statistic and $Z_t^{(K)}$ is the smallest LLR at time $t$. 
 
 When  the size of the affected subset is known in advance, i.e., $\um=\om=m$,  we have
\begin{equation} \label{qq1}
 \hat{Z}_t=  \sum_{k=1}^{m}  Z_{t}^{(k)} .
 \end{equation}
 Indeed, for any  $\cA  \in \cP_{m}$ we have 
$
Z_t^\cA \leq \sum_{k=1}^{m}  Z_{t}^{(k)}.
$
Therefore,
$
\hat{Z}_t \leq   \sum_{k=1}^{m}  Z_{t}^{(k)},
$
and the  upper bound is attained by the subset which consists of  the $m$ channels with the highest LLR values at time $t$. 

In the more general case that $\um<\om$ we have 
$$
\hat{Z}_t= \sum_{k=1}^{\um}  Z_{t}^{(k)} + \sum_{k=\um+1}^{\om} ( Z_{t}^{(k)} )^{+},
 $$
and the \mbox{G-SLRT} takes the following form:
 \begin{align} \label{gsprt_ordered}
\begin{split}
 \hat{\tau}&= \inf\set{ t\geq 0: \sum_{k=1}^{\om}  \left(  Z_{t}^{(k)}  \right)^{+} \geq b \quad  \text{or}  \quad
 \sum_{k=1}^{\um}  Z_{t}^{(k)}   \leq -a } \\
\hat{d}&=  \begin{cases}
1 & \; \text{when}  \quad \sum_{k=1}^{\om}  ( Z_{\hat{\tau}}^{(k)} )^{+}\geq b  \\	
0  & \; \text{when}  \quad  \sum_{k=1}^{\um}  Z_{\hat{\tau}}^{(k)}  \leq -a  \\
\end{cases}.
\end{split}
\end{align}
Indeed,   for any $\cA \in \cP_{\um, \om} $ we have 
$$
Z^{\cA}_t \leq \sum_{k=1}^{\um}   Z_{t}^{(k)}   +   \sum_{k=\um+1}^{\om}  (Z_{t}^{(k)} )^{+},
$$
and the upper bound is attained by the subset  which consists of the $\um$ channels with the top  $\um$ LLRs and the next (if any)  top  $\om-\um$ channels that have positive LLRs.   

\subsection{Generalization}\label{ssec:Generalization}

It is possible to generalize the GLR detection statistic in \eqref{gslrt} by applying different weights to the  LLRs of the various hypotheses. Specifically, let  $\cP$ be an arbitrary class and $\{p_{\cA}\}_{\cA \in \cP}$ an arbitrary family of positive numbers  (weights) that add up to 1.  Then, the weighted GLR detection statistic may be defined as 
\begin{equation} \label{gslrt3}
  \max_{\cA \in \cP}  \left(  Z^{\cA}_t + \log p_{\cA}  \right) .
\end{equation}
It is straightforward to see that the asymptotic optimality properties that we established in the previous section remain valid for any selection of weights (that do not depend on the 
thresholds or the error probabilities).  Moreover, the resulting sequential test is as  feasible as the \mbox{G-SLRT}, as long as there are positive numbers  $\{ p_{k} \}_{1 \leq k \leq K}$ such that each $p_\cA$ is 
proportional to  $\prod_{k \in \cA} p_k$, i.e., 
\begin{equation} \label{weights}
p_\cA = C(\cP) \, \prod_{k \in \cA} p_k, \quad 
C(\cP) = \left(\sum_{\cA \in \cP} \prod_{k \in \cA} p_k \right)^{-1},
 \end{equation}
that is, $C(\cP)$ is a normalizing constant.  Indeed, in this case, the weighted GLR statistic \eqref{gslrt3}  takes the form
\begin{equation*} 
 \max_{\cA \in \cP}  \sum_{k \in \cA} \left( Z^{k}_t + \log p_{k}   \right)+ \log C(\cP)
 \end{equation*}
and the discussion in Subsection \ref{ssec:GSLRT} applies with $Z^{k}_t$ replaced by $Z^{k}_t+\log p_k$ and thresholds $a$ and $b$ replaced by
 $a+ \log C(\cP)$ and $b- \log C(\cP)$, respectively.

\section{Mixture-based Sequential Likelihood Ratio Test} \label{s:WSLRT}

In this section, we propose an alternative sequential test that is based on 
averaging, instead of maximizing,  the likelihood ratios that correspond to the  different hypotheses. We show that it has the same asymptotic optimality properties and similar feasibility as the \mbox{G-SLRT}.

\subsection{Definition and Error Control}\label{ssec:WSLRTdefin}

Let  $\cP$ be an arbitrary class, $\{p_{\cA}\}_{\cA \in \cP}$ an arbitrary family of positive numbers that add up to 1 (weights) and  consider the probability  measure 
 \begin{equation} \label{mixture_measure}
 \overline{\Pro} := \sum_{\cA \in \cP} p_{\cA} \Pro^{\cA}.
\end{equation}
Then, the Radon-Nikod\'ym derivative of  $\overline{\Pro}$ versus $\Pro_0$ given $\cF_t$ is 
\begin{align} \label{mixture}
 \overline{\Lambda}_{t} &:= \frac{d \overline{\Pro}}{ d\Pro_0} \Big|_{\cF_t} = \sum_{\cA \in \cP} p_{\cA} \Lambda^{\cA}_t  =  \sum_{n=1}^{K} \sum_{\cA \in \cP \cap \cP_{n}}   p_{\cA} \Lambda^{\cA}_t.
\end{align}
If we replace the generalized likelihood ratio statistic $\hat{Z}$ in \eqref{gslrt} by  the logarithm of the  mixture likelihood ratio,  $\overline{Z}_t := \log \overline{\Lambda}_{t}$, then we obtain the following sequential test: 
\begin{align} \label{wslrt}
\begin{split}
 \overline{\tau} &= \inf\set{ t: \overline{Z}_{t}  \notin (-a,b)}, \quad 
 \overline{d} :=  \begin{cases}
1 & \; \text{when}  \quad \overline{Z}_{\overline{\tau}} \geq b   \\	
0  & \; \text{when}  \quad  \overline{Z}_{\overline{\tau}} \leq -a \\
\end{cases} ,
\end{split}
\end{align}
to which we refer as the Mixture Sequential Likelihood Ratio Test (\mbox{M-SLRT}). In the following lemma we show how to select the thresholds in order to guarantee the desired error control for \mbox{M-SLRT}.

\begin{lemma}\label{Lem:PEuppermlrt}
For any  positive thresholds $a$ and $b$ we have 
\begin{equation}\label{PEupper_mix}
 \Pro_0(\overline{d} =1) \le  e^{-b} \quad \text{and} \quad \max_{\cA\in \cP} \Pro^\cA(\overline{d}=0) \le \left( \min_{\cA \in \cP} p_{\cA} \right)^{-1}  e^{-a}.
\end{equation}
Therefore,  for any $\alpha, \beta \in (0,1)$, $(\overline\tau,\overline{d}) \in \ccab(\cP)$ when the thresholds are selected as follows:
\begin{equation} \label{selection_mix}
 b=|\log \alpha| \quad \text{and} \quad a=|\log \beta| -\min_{\cA \in \cP} (\log p_{\cA}).
 \end{equation}
\end{lemma}

\begin{IEEEproof} 
Let $\overline{\Exp}$ be the  expectation that corresponds to  the mixture measure   $\overline{\Pro}$ defined in \eqref{mixture_measure}.  Since $\overline{Z}_{\overline\tau} \ge b$ on $\{\overline{d} =1\}$, from Wald's likelihood ratio identity we have 
\[
\begin{aligned} 
\Pro_0(\overline{d} =1)   &= \overline{\Exp} \brcs{\exp\{-\overline{Z}_{\overline{\tau}}\} ;  \overline{d} =1}  \le  e^{-b} ,
\end{aligned}
\]
which proves the first  inequality in \eqref{PEupper_mix}. In order to prove the second inequality we note that, for  any $\cA \in \cP$,  on the event  $\{\overline{d} =0\}$  we have  $-a\geq \overline{Z}_{\overline{\tau}}\geq  Z^\cA_{\overline{\tau}} + \log p_\cA$. Consequently, from Wald's likelihood ratio identity we obtain
\begin{align*}
\Pro^\cA(\overline{d}=0) &= \Exp_0\brcs{ \exp \{ Z^\cA_{\overline{\tau} } \} ; \;  \overline{d} =0 }  \leq   p_{\cA}^{-1} e^{-a} .
\end{align*}
Since this inequality is true for any $\cA \in \cP$,  maximizing both sides with respect to $\cA$ proves the second  inequality in \eqref{PEupper_mix}. 
\end{IEEEproof} 

\subsection{Asymptotic Optimality}\label{ssec:AOMSLRT}

The following theorem shows that the \mbox{M-SLRT}   has exactly the same  asymptotic optimality properties as the \mbox{G-SLRT}.

\begin{theorem}\label{Th:AOCI}
Consider an arbitrary class of possibly affected subsets, $\cP$, and 
suppose that the  thresholds of the \mbox{M-SLRT} are   selected according to \eqref{selection_mix}.  If   $r$-complete convergence conditions \eqref{rcompleteLLRs} hold, then 
for all $1 \leq q\le r$ we have as $\alpha_{\max}\to 0$:
 \begin{align}
\Exp_{0}[\overline{\tau}^q] & \sim \brc{\frac{|\log \beta|}{\min\limits_{\cA \in \cP}   I_{0}^{\cA} } }^q \sim \inf_{(\tau,d) \in \ccab(\cP)} \Exp_{0}[\tau^q] , \label{new_sprt1}
\\
 \Exp^{\cA}[\overline{\tau}^q] &\sim \brc{\frac{|\log \alpha|}{I_{1}^{\cA} }}^q \sim \inf_{(\tau,d) \in \ccab(\cP)} \Exp^{\cA}[\tau^q]  \quad \text{for every}~ \cA \in \cP. \label{new_sprt2}
 \end{align}
 Moreover, if the LLRs $Z_t^k$ have independent increments,   then the asymptotic relationships 
 \eqref{new_sprt1}--\eqref{new_sprt2}  hold for every $q > 0$   as long as  the  almost sure convergence conditions \eqref{as} are satisfied. 
\end{theorem} 
 
 \begin{IEEEproof}
 The proof is based on the observation that for every $t \in \bN$ we have 
 \begin{align} \label{pathwise}
 \min_{\cA \in \cP} (\log p_{\cA}) &\leq \overline{Z_t} - \hat{Z}_t \leq     \max_{\cA \in \cP} (\log p_{\cA}) +\log |\cP|.
  \end{align}
  \end{IEEEproof}
 
 \subsection{Feasibility }

Similarly to the \mbox{G-SLRT}, the \mbox{M-SLRT} is computationally feasible even when $K$ is large if the weights are selected according to \eqref{weights}.  Then, the mixture likelihood ratio takes the form
\begin{align*}
\overline{\Lambda}_{t}  &=  C(\cP)   \, \sum_{m=1}^{K} \sum_{\cA \in \cP \cap \cP_{m}} \, \prod_{k \in \cA}  \left( p_k \Lambda_t^{k} \right)  .
\end{align*} 
When in particular there  is an upper and a lower bound on the size of the affected subset, i.e.,  $\cP=\cP_{\um, \om}$  for some $1 \leq \um \leq \om \leq K$,    the mixture likelihood ratio statistic  takes the form
\begin{align} \label{mixture2}
\overline{\Lambda}_{t} =    C(\cP) \,  \sum_{m=\um}^{\om}  \sum_{\cA \in \cP_{m}} \prod_{k \in \cA}  \left( p_k \Lambda_t^{k} \right) 
\end{align} 
and its computational  complexity  is  polynomial in the number of channels, $K$.  However, in the special case of complete uncertainty ($\um=1, \om=K$), the \mbox{M-SLRT} requires only $\calo(K)$ operations. Indeed, if we set for simplicity   $p_k=p$ and $ \pi=p/(1+p)$,   then  the mixture likelihood ratio in \eqref{mixture2}  admits the following representation  for the class $\cP=\overline{\cP}_{K}$:
\begin{equation} \label{hatid}
  \overline{\Lambda}_{t}=    C(\cP)  \,  [ (1-\pi)^{-K} \tilde{\Lambda}_{t}-1]
  \end{equation}
 where  the statistic $\tilde{\Lambda}_t$ is defined as follows:
\begin{equation} \label{hat}
\tilde{\Lambda}_{t}:=
 \prod_{k=1}^{K}  \left( 1-\pi+\pi\, \,  \Lambda_{t}^{k} \right).
\end{equation}

\begin{remark}
The statistic $\tilde{\Lambda}_{t}$ has an appealing statistical interpretation, as it is  the   likelihood ratio that corresponds to the case that  each  channel belongs to the  affected subset with  probability $\pi \in (0,1)$.  It is possible to use $\tilde{\Lambda}_{t}$ as the detection statistic and incorporate prior  information by an appropriate selection of $\pi$. For instance, if we know  the exact size of the affected 
subset, say $\cP=\cP_m$, we may set    $\pi =m/K$, whereas if we know that  at most $m$ channels may be affected, i.e.,  $\cP=\overline{\cP}_m$, then we may set  $\pi =m/(2K)$. 
This approach was consider in \cite{felsok,Xie&Siegmund-AS13} for a multistream quickest change detection problem.
\end{remark}

\section{Examples} \label{sec:Examples}

In this section, we consider three particular examples to which the previous results apply.

\subsection{ A Linear Gaussian State-Space Model}\label{ssec:Kalman}

First, we present the  example of a linear state-space (hidden Markov) model,   in which the LLR process, $\{Z_t^{k}\}$, has independent increments and  Theorems \ref{Th:GLRTAOindep_gslrt} and \ref{Th:AOCI} are applicable. Let $X_t^k=(X_{t,1}^k,\dots, X_{t,\ell}^k)^\top$ be the  $\ell$-dimensional observed vector  in the $k$-th channel  at time $t$ and let 
$\theta_t^k=(\theta_{t,1}^k, \dots,\theta_{t,m}^k)^\top$  be the unobserved $m$-dimensional Markov vector and suppose that
\[
\begin{aligned}
\theta_t^k & = F^k \, \theta_{t-1}^k + W_{t-1}^k + i \, b_\theta^k, \quad  \theta_0^k=0,
\\
X_t^k & = H^k \, \theta_t^k + V_t^k + i \, b_x^k ,
\end{aligned}
\]
where $W_t^k$ and~$V_t^k$ are zero-mean Gaussian i.i.d.\ vectors having covariance matrices  $K_W^k$ and~$K_V^k$, respectively; $b_\theta^k=(b_{\theta,1}^k,\dots,b_{\theta,m}^k)^\top$ 
and $b_x^k=(b_{x,1}^k,\dots,b_{x,\ell}^k)^\top$ are the mean values; $F^k$ is the $(m\times m)$ state transition matrix; $H^k$ is the $(\ell \times m)$ matrix, 
and the index $i=0$ if the mean values in the $k$-th channel (component) are not affected and $i=1$ otherwise.

It can be shown that under the null hypothesis $\Hyp_0^k$  the observed sequence~ $X^k$ has an equivalent representation 
\[
X_t^k = H^k \hat{\theta}_t^k + \xi_t^k, \quad t\in\Nbb
\]
with respect to the ``innovation'' sequence  $\xi_t^k=X_t^k - H^k \hat{\theta}_t^k$,  where $\xi_t^k \sim \cN(0,\Sigma_t^k)$, $t =1,2,\ldots\,$ are independent Gaussian vectors and 
$\hat{\theta}_t^k=\Exp_0^k[\theta_t^k\,|\, X_1^k, \dots, X_{t-1}^k]$ 
 is the optimal one-step ahead predictor in the mean-square sense, i.e., the estimate of~$\theta_t^k$ based on observing $X_1^k,\dots,X_{t-1}^k$, which can be obtained by the Kalman filter 
 (cf., e.g., \cite{balakrishnan-book87}).  On the other hand, 
 under $\Hyp_1^k$ the observed sequence $X^k$ admits the following
 representation 
$$
X_t^k = \Upsilon_t^k+ H^k \hat{\theta}_t^k + \xi_t^k, \quad t \in \Nbb,
$$
where $\Upsilon_t^k$ depends on~$t$ and can be computed using relations given in~\cite[pp.~282-283]{Basseville&Nikiforov-book93}. Consequently, 
 the local LLR $Z^k$  can be written as
$$
Z_t^k  = \sum_{s=1}^t (\Upsilon_s^k)^\top (\Sigma_s^k)^{-1} \xi_s^k - \frac{1}{2}\sum_{s=1}^t (\Upsilon_s^k)^\top (\Sigma_s^k)^{-1} \Upsilon_s^k,
\quad t \in \bN,
$$
where $\Sigma_t$, $t \in \bN$  are given by Kalman's equations (see \cite[Eq.~(3.2.20)]{Basseville&Nikiforov-book93}). Thus, each $Z^k$
has independent Gaussian increments. Moreover, it is easily seen that the normalized LLR~$t^{-1} Z_t^k$ converges almost surely   as $t \to\infty$ to $I_1^k$ under $\Pro_1^k$ and  $-I_1^k$ under $\Pro_0^k$, where 
$$
I_1^k= \frac{1}{2} ~\lim_{t\to\infty} ~\frac{1}{t}\sum_{s=1}^t (\Upsilon_s^k)^\top (\Sigma_s^k)^{-1} \Upsilon_s^k .
$$ 
Therefore, by  Theorem~\ref{Th:GLRTAOindep_gslrt} 
and Theorem~\ref{Th:AOCI}, the \mbox{G-SLRT} and the \mbox{M-SLRT} are asymptotically optimal  with respect to all moments of the sample size. 

\subsection{An Autoregression Model with Unknown Correlation Coefficient}\label{ssec:AR}

Suppose that the observations in the channels are Markov Gaussian (AR(1)) processes  of the form
\[
X_t^k = \rho^k X_{t-1}^k + \xi_t^k, \quad t \in \bN, \quad X_0^k=0,
\] 
where $\{\xi_t^k\}_{t\in \Nbb}$, $k=1,\dots,K$ are mutually independent sequences of  i.i.d.\ normal random variables with  zero mean and unit variance. Suppose that  $\rho^k=\rho_i^k$ under $\Hyp_i^k$, $i=0,1$, where 
$\rho_i^k$  are known constants. Then, the transition densities are $f_i^k(X_t^k|X_{t-1}^k) = \varphi(X_t^k-\rho_i^kX_{t-1}^k)$, $i=0,1$,  where  $\varphi$  is the density of the standard normal distribution, and the LLR in the $k^{th}$ channel can be written as $Z_t^k = \sum_{s=1}^t g_k(X_s^k,X_{s-1}^k)$, where 
\begin{equation}\label{g}
g_k(y,x) :=  \log \left( \frac{\varphi(y-\rho_1^k x)}{\varphi(y-\rho_0^k x}  \right) =  \frac{1}{2} \left[(y-\rho_0^k x)^2 -(y-\rho_1^k x)^2 \right]=
(\rho_1^k-\rho_0^k)x \left[ y-\frac{\rho_1^k+\rho_0^k}{2}  x \right].
\end{equation}
In order to show that  $\{t^{-1} Z_t^k\}$ converges asymptotically as $t \rightarrow \infty$, let us further assume that  $|\rho_i^k| <1$, $1 \leq k \leq K$, $i=0,1$,  so that $X^k$ is stable. Let $\lambda_i^k$ be  the invariant distribution of $X^k$ under $\Hyp_i^k$, which  coincides with the  distribution of 
\begin{equation}\label{sec:Ex.1-1n}
w_i^k=\sum^{\infty}_{t=1} (\rho_i^k)^{t-1}\,\xi_t^k, \quad i=0,1.
\end{equation}
By  a slight extension of Theorem~5.1 in \cite{PergTar} to $r>1$ (see Appendix \ref{appen;AR}),  it can be shown that  under $\Pro_1^k$ the normalized LLR process $\{t^{-1} Z_t^k\}$ converges  as $t\to\infty$  $r$-completely  for every  $r \ge 1$ to 
$$
I_1^k=  \int_{-\infty}^\infty \brc{\int_{-\infty}^\infty g_k(y, x) \,\varphi(y-\rho_1^k x)d y}\, \lambda_1(d x).
$$
 In the Gaussian case  considered, $\lambda_1^k$ is  $\cN(0,(1-\rho_1^k)^{-2})$, so $I_1^k$  can be calculated explicitly as
$I_1^k = (\rho_1^k - \rho_0^k)^2/2[1-(\rho_1^k)^2]$. By symmetry, 
under $\Pro_0^k$ the normalized LLR $\{t^{-1} Z_t^k\}$ converges $r$-completely for all $r \ge 1$ to~$-I_0^k$ with  $I_0^k = (\rho_1^k - \rho_0^k)^2/2[1-(\rho_0^k)^2]$. Thus, by Theorem~\ref{Th:GLRTAO} and Theorem~\ref{Th:AOCI}, both tests, the \mbox{G-SLRT} and the \mbox{M-SLRT}, are asymptotically optimal minimizing all moments of the stopping time distribution.

\subsection{Multichannel Invariant Sequential $t$-Tests}\label{ssec:ttest}

Suppose that the observations in channels have the form
\[
X_t^k = i \mu_k + \xi_t^k, \quad t \in \bN, ~~ 1 \leq k \leq K,
\]
where  $\xi_t^k\sim \cN(0,\sigma_k^2)$, $t\in \bN$ are  zero-mean, normal i.i.d.\ (mutually independent)  sequences (noises)  with unknown variances $\sigma_k^2$.   Under the local null hypothesis in the $k^{th}$ stream, $\Hyp_0^k$, there is no signal in the $k^{th}$ stream ($i=0$). Under the local alternative hypothesis in the $k^{th}$ stream, there is a signal $\mu_k>0$ in the $k^{\rm th}$ channel. Therefore, the hypotheses $\Hyp_0^k$, $\Hyp_1^k$ are not simple and our results cannot be directly applied. Nevertheless,  if we assume that the 
value of the ``signal-to-noise'' ratio $Q_k=\mu_k/\sigma_k$ is known, we can transform this into a testing problem of simple hypotheses in the channels by  using the principle of invariance, since the problem  is invariant under the group of scale changes. Indeed, the maximal invariant statistic in the $k^{\rm th}$ channel is 
$\bY_t^k=(1,  X_2^k/X_1^k, \dots,  X_t^k/X_1^k)$ and it can be shown \cite[Sec 3.6.2]{TNB_book2014} that the invariant LLR, which is built based on the maximal invariant $\bY_t^k$, is given by $Z_t^k =  \log [J_t(Q_k T_t^k)/J_t(0)]$,  where 
\begin{equation}\label{t_statistic}
T_t^k  = \frac{t^{-1}\sum_{s=1}^t X_j^k}{\left\{t^{-1}\sum_{s=1}^t (X_s^k)^2\right \}^{1/2}} 
\end{equation}
and 
\begin{equation*}
J_t(z) =\int_0^{\infty} \frac{1}{u}\exp\set{\brcs{-\frac{1}{2}u^2 + z u + \log u} \, t }\,d u.
\end{equation*}

Note that $T_t^k$ is the Student $t$-statistic, which is the basis for Student's $t$-test in the fixed sample size setting.  For this reason, we refer to the sequential tests \eqref{gslrt} and \eqref{wslrt} that are based on the invariant LLRs as 
$t$-tests, in particular as the \mbox{$t$-G-SLRT} and the \mbox{$t$-M-SLRT}, respectively. Although the  invariant
LLR~$Z_t^k$ is difficult to calculate explicitly, it  can be approximated by  $g_k(T_t^k) \, t$,  using a uniform version of the Laplace asymptotic integration technique,  where the function $g_k(x)$ is given by 
\[
g_k(x)  = \frac{1}{4} x \brc{x+\sqrt{4+x^2}}+\log\brc{x+\sqrt{4+x^2}} - \log 2 - \frac{1}{2} Q_k^2 , \quad x \in \bR.
\]
Indeed, as shown in \cite[Sec 3.6.2]{TNB_book2014},  there is a  finite positive constant $C$ such that for all $t\ge 1$ we have $|Z_t^k-g_k(T_t^k) \, t|\le C$, or equivalently,
\begin{equation}\label{LLRapproxinv}
\abs{t^{-1} Z_t^k -g_k(T_t^k)}  \leq C/t.
\end{equation}
It follows from \eqref{LLRapproxinv} that if  under~$\Pro_i^k$ the $t$-statistic~$T_t^k$ converges $r$-completely   to a constant~$V_i^k$ as $t \rightarrow \infty$, then the normalized LLR $t^{-1} Z_t^k$ converges in a similar sense to~$g_k (V_i^k)$, $i=0,1$. Therefore,  it suffices to study the limiting behavior of~$T_t^k$.  Since  for every $r\ge 1$ we have  $\Exp_i^k[\abs{X_1^k}^r]<\infty$, $i=0,1$, for every $r\ge 1$ we obtain 
\begin{equation*}
T_t^k \xrightarrow[t\to\infty]{\Pro_1^k-r-\text{completely}} \frac{\Exp_1^k [X_1^k]}{\sqrt{\Exp_1^k [(X_1^k)^2]}} = \frac{Q_k}{\sqrt{1+Q_k^2}}, \qquad T_t^k \xrightarrow[t\to\infty]{\Pro_0^k-r-\text{completely}} 0,
\end{equation*}
which implies that  the $r$-complete convergence condition~\eqref{rcompleteLLRs} for the  normalized LLR $\{t^{-1} Z_t^k\}$  holds  for all $r\ge 1$ with
\begin{equation*}
I_1^k = g_k\brc{\frac{Q_k}{\sqrt{1+Q_k^2}}} \quad \text{and} \quad I_0^k = \frac{1}{2} Q_k^2 .
\end{equation*}
It is easy to verify that $I_1^k>0$ and $I_0^k>0$.  Hence, by Theorem~\ref{Th:GLRTAO} and Theorem~\ref{Th:AOCI}, the invariant \mbox{$t$-G-SLRT}  and \mbox{$t$-M-SLRT}  asymptotically minimize all moments of the stopping time distribution.

\section{Simulation Experiments} \label{sec:Simulations}

In this section  we  present the results of a simulation study whose goal is to compare the  performance of  the \mbox{G-SLRT} and  the \mbox{M-SLRT}, as well as to quantify the effect of prior information on the detection performance. 

\subsection{ Computation of Error Probabilities Via Importance Sampling}\label{ssec:EPIS}

Since the  type-I and type-II errors for both the \mbox{G-SLRT} and the \mbox{M-SLRT} correspond to ``rare events'',    we  rely on importance sampling for the  computation of  these  probabilities. We illustrate this method for the  \mbox{G-SLRT}, since the approach for  the \mbox{M-SLRT} is identical. 
 
 We start with the maximal type-II error. From \eqref{com} it follows that
 for every $\cA \in \cP$ we have
$$
 \Pro^\cA(\hat{d}=0) =   \Exp_0\brcs{\exp\{Z^\cA_{\hat\tau}\} ; \hat{d} =0}.
$$
Therefore,   all probabilities $\Pro^\cA(\hat{d}=0)$,  $\cA \in \cP$, can be computed simultaneously by  simulating the \mbox{G-SLRT}, $(\hat{\tau}, \hat{d})$,  under $\Pro_0$, which then allows the computation of the maximal type-II error probability.   This computation is particularly simplified when all hypotheses are identical, in the sense that $\Pro_i^k$ does not depend on $k$, $i=0,1$. Indeed, in this case, 
$$
\max_{\cA \in \cP} \Pro^\cA(\hat{d}=0) =  \max_{1 \leq m \leq K:  \cP \cap \cP_{m} \neq \varnothing} \,  \Pro^{\cA_{m}} (\hat{d}=0),
$$
where  $\cA_m$ is an arbitrary set in $\cP_m$, say $\cA_m= \{1, \ldots,m\}$. When in particular $\cP=\cP_{\um, \om}$ for some $1 \leq \um \leq \om\leq K$, then 
$$
\max_{\cA \in \cP} \Pro^\cA(\hat{d}=0) = \max_{\um \leq m \leq \om } \Pro^{\cA_m}  (\hat{d}=0).
$$

We now turn to the  computation of the type-I error probability for which  we rely on the change of measure  $\Pro_0 \rightarrow \overline{\Pro}$, where $\overline{\Pro}$ is the mixture probability measure defined in 
 \eqref{mixture_measure} with uniform weights, i.e., $p_{\cA}=(\log |\cP|)^{-1}$ for every $\cA \in \cP$.   Indeed,  from Wald's likelihood ratio identity it follows that 
\begin{equation} \label{is}
\Pro_0(\hat{d}=1) =   \overline{\Exp} \brcs{\overline{\Lambda}^{-1}_{\hat\tau} ; \hat{d} =0}  =   \overline{\Exp} \brcs{\exp\{-\overline{Z}_{\hat\tau}\} ; \hat{d} =0} ,
\end{equation}
where $\overline{\Exp}$ refers to expectation under $\overline{\Pro}$. Even though the test statistic does not coincide with the likelihood ratio statistic that is used for the change of measure,  the second moment (and, consequently, the variance)  of this estimator  is bounded above by
$$ 
\overline{\Exp} \brcs{\exp\{-2\overline{Z}_{\hat\tau} \}; \hat{d} =0}  \leq 
\overline{\Exp} \brcs{\exp\{-2(\overline{Z}_{\hat\tau} -\hat{Z}_{\hat\tau}) -2\hat{Z}_{\hat\tau}  \} ; \hat{d} =0}  \leq \exp\{2 \log |\cP| -2b \}  ,
$$
since from \eqref{pathwise} it follows that   $\hat{Z}_t -  \overline{Z}_t\leq \log |\cP|$ for every $t$. Lemma \ref{Lem:PEupperglrt} implies that  $\Pro_0(\hat{d}=1) \leq  \exp\{-b + \log |\cP| \} $ for every $b$. If also there is some constant $c \in (0,1)$ such that  $\Pro_0(\hat{d}=1) \sim c  \exp\{-b + \log |\cP| \} $ as $b \rightarrow \infty$, then the  relative error of this importance sampling estimation is asymptotically bounded as $b \rightarrow \infty$ since
$$\sqrt{\overline{\Exp} \brcs{\exp\{-2\overline{Z}_{\hat\tau} \}; \hat{d} =0}}  =\calo(\Pro_0(\hat{d}=1)).
$$

As far as the computational complexity of this computation concerns, from  the definition of  $\overline{\Pro}$ we have  that the expectation in \eqref{is} can be written as follows:  
$$
 \sum_{\cA \in \cP} p_{\cA} \, \Exp^{\cA} \brcs{ 
 \frac{1} {\overline{\Lambda}_{\hat\tau}} ; \hat{d} =0} 
= \sum_{m=1}^{K}  \sum_{\cA \in \cP_{m} \cap \cP} p_{\cA} \, \Exp^{\cA} \brcs{  \frac{1} {\overline{\Lambda}_{\hat\tau}} ; \hat{d} =0} ,
$$
which  requires simulating the \mbox{G-SLRT} under each $\Pro^{\cA}$ with $\cA \in \cP$. This computation is considerably simplified in the case of symmetric hypotheses,   in which case the expectation in  \eqref{is} takes the form:
$$
 \sum_{m=1}^{K}  \frac{|\cP_{m} \cap \cP|}{|\cP|}  \, 
  \Exp^{ \cA_m } \brcs{\Lambda^{-1}_{\hat\tau}; \hat{d} =0},
$$
where  $\cA_m= \{1, \ldots,m\}$. When in particular we have a class of the form $\cP_{\um, \om}$, then the expectation in \eqref{is} becomes
$$
  \sum_{m=\um}^{\om}  \frac{|\cP_{m}|}{|\cP|}  \, 
  \Exp^{ \cA_m } \brcs{\Lambda^{-1}_{\hat\tau} ; \hat{d} =0},
$$
which requires simulating the \mbox{G-SLRT}  under only $\om-\um$ scenarios.

\subsection{A Simulation Study for an Autoregressive Model}\label{ssec:SimulationsAR}

We now present the results of a simulation study in the context of the  autoregression of Subsection~\ref{ssec:AR}. We assume that the hypotheses are symmetric in the sense that   $\rho_0^k=0$ and $\rho_1^k=\rho=0.5$, therefore the Kullback-Leibler divergences take the form $I_1^k=I_1= (1/2) \rho^2/(1-\rho^2)$, $I_0^k=I_0= (1/2) \rho^2$.  Our goal is to compare the \mbox{G-SLRT} against the \mbox{M-SLRT} (with uniform weights) for two different scenarios regarding the available prior information; in the first one, the size of the  affected subset is assumed to be known, i.e.,  $\cP=\cP_m$ where  $m$ is  the cardinality of the true affected subset; in the second one,  there  is  complete uncertainty regarding the affected subset ($\cP=\overline{\cP}_K$). 

We assume that  $\alpha=\beta$, where  $\alpha, \beta$ are the desired type-I and type-II error probabilities for the two tests. For the \mbox{G-SLRT} we select the pair of  thresholds, $a, b$ such  that $b=a+ \log |\cP|$,  where $\cP$ is the class of possibly affected subsets. Then, from  \eqref{selection}  it follows that both error probabilities will be bounded above by $\exp(-a)$.  In Tables  \ref{tab1} and \ref{tab2} we present the operating characteristics of the \mbox{G-SLRT} when $a= 8.2$, in which case  both error probabilities are bounded  by $\exp(-a)=2.75\cdot 10^{-4}$. 
\begin{table}[h]   
\centering
\caption{Error probabilities and expected sample under the null of 
the  \mbox{G-SLRT} and  the \mbox{M-SLRT}. The \mbox{G-SLRT}
 thresholds are $a=8.2$ and   $b=a+ \log |\cP|$, whereas  the  \mbox{M-SLRT} thresholds are $b=8.2$ and  $a=b+ \log |\cP|$, where $\cP$ is the class of possibly affected subsets. 
Standard errors are presented in parentheses based on 1,600 simulation runs.  \label{tab1}  } \vspace{-3mm} 
\begin{tabular}{|c|c|c|c|c|c|c|}
\hline
$\cP$ & \multicolumn{2}{c|}{ $\Pro_0(d=1)$} & \multicolumn{2}{c|}{ $\Exp_0[T]$} & \multicolumn{2}{c|}{ $\max_{\cA \in \cP} \Pro^\cA(d=0)$} \\ 
& G-SLRT  & M-SLRT  &G-SLRT  & M-SLRT   &G-SLRT  & M-SLRT \\ \hline \hline
 $\overline{\cP}_{K}$ &   2.39 (0.023)   $\cdot10^{-5}$ & 9.3 (0.06)   $\cdot10^{-5}$&  140.6  (0.9) &	 146.1 (0.9) & 2.12 (0.12)  $\cdot10^{-5}$  & 1.40 (0.09)  $\cdot 10^{-5}$\\
 $\cP_{1}$   & 3.33 (0.08) $\cdot 10^{-5}$  &1.21 (0.01) $\cdot 10^{-4}$ 
 &  140.0 (0.86) & 120.0 ( 0.8) & 2.15 (0.12)
$\cdot 10^{-5}$& 1.56 (0.09) $\cdot 10^{-4}$		\\
                $\cP_{3}$   & 7.17  (0.12)  $\cdot 10^{-5}$   &1.02 (0.014) $\cdot 10^{-4}$ &	54.6 (0.4) 	&	41.5  (0.3) 	& 3.90  (0.48) $\cdot 10^{-6}$	&	1.03 (0.11)  $\cdot 10^{-4}$ \\
                $\cP_{6}$   &  2.78 (0.07)  $\cdot 10^{-5}$&  8.9 (0.14) $\cdot 10^{-5}$ & 24.4  (0.2)  &	 17.8 (0.2) & 2.67 (0.31)                $\cdot 10^{-6}$    & 1.02 (0.15) $\cdot 10^{-4}$  \\
                $\cP_{9}$   &  3.33 (0.085)  $\cdot 10^{-5}$& 7.38 (0.15)    $\cdot 10^{-5}$ &  11.4 (0.1) 		&	9.6 (0.1) & 1.96 (0.08) $\cdot 10^{-5}$  & 9.09 (0.48) $\cdot 10^{-5}$	 \\ \hline
\end{tabular}
\end{table}

\begin{table} [h] 
\centering
\caption{Same setup as in Table \ref{tab1}. Here, we report the expected sample size of each test under the alternative hypothesis for various scenarios regarding the  number of affected channels \label{tab2}}
\vspace{-3mm}
\begin{tabular}{|c|c|c|c|c|}
\hline
  $|\cA|$  & \multicolumn{4}{c|}{ $\Exp^\cA[T]$  }  \\ 
 & \multicolumn{2}{c|}{ G-SLRT } &  \multicolumn{2}{c|}{ M-SLRT }  \\
  &  $\overline{\cP}_{K}$  & $\cP_{|\cA|}$  &$\overline{\cP}_{K}$  & $\cP_{|\cA|}$\\
  \hline \hline
  1 &   100.0 (1.2) &  74.5  (1.0)	&  96.0 (1.2) &   71.2 (0.9) \\
  3 & 33.9 (0.4)  & 29.8 (0.4) &  29.9 (0.4) &	 28.4 (0.4)	\\
  6 & 	17.25 (0.18) &  15.7 (0.2) &	 14.82 (0.18) 	&	14.3 (0.2)	\\
  9 & 12.0 (0.1) & 9.7 (0.1)   &  9.87 (0.10)  & 8.9 (0.1) \\ \hline
\end{tabular}
\end{table}

We select the pair of thresholds $a,b$ of the \mbox{M-SLRT} such that  $a=b+ \log |\cP|$. Then, from \eqref{selection_mix}  it follows that  both error probabilities will be bounded above by $\exp(-b)$.  In Tables \ref{tab1} and \ref{tab2} we present the operating characteristics of the \mbox{M-SLRT} when $b= 8.2$, in which case  both error probabilities are also bounded by  $\exp(-b)=2.75 \cdot 10^{-4}$. In this way, the results for the two schemes  are comparable.  From these tables we can see that  in all cases the actual error probabilities  are much smaller than the target value of  2.75 $\cdot 10^{-4}$,  but this upper bound is much more conservative for the G-SLRT than for the M-SLRT.

\begin{figure}
\centering
\begin{tabular}{cc}
\includegraphics[width=0.5\linewidth, height=0.5\linewidth]{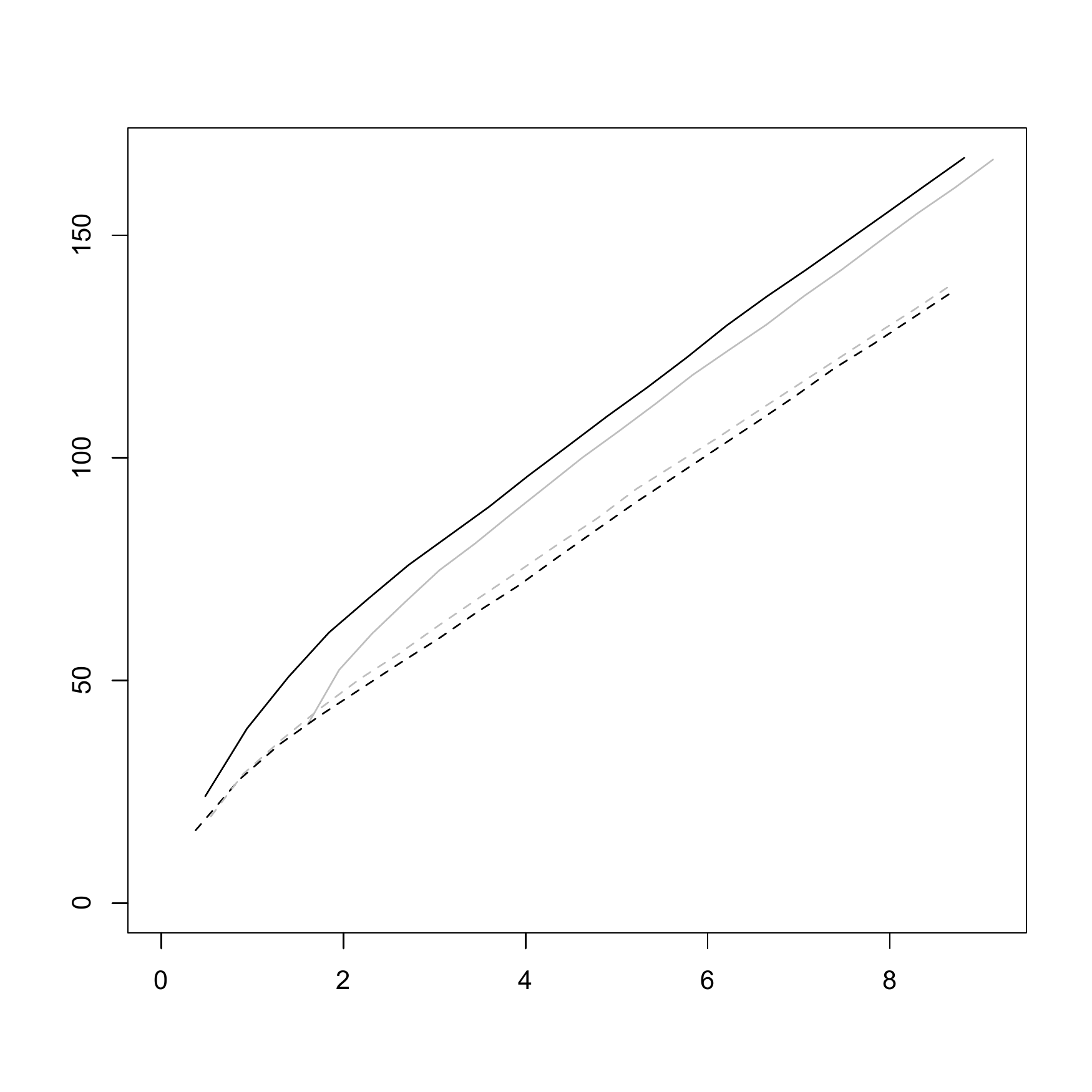}
&
\includegraphics[width=0.5\linewidth, height=0.5\linewidth]{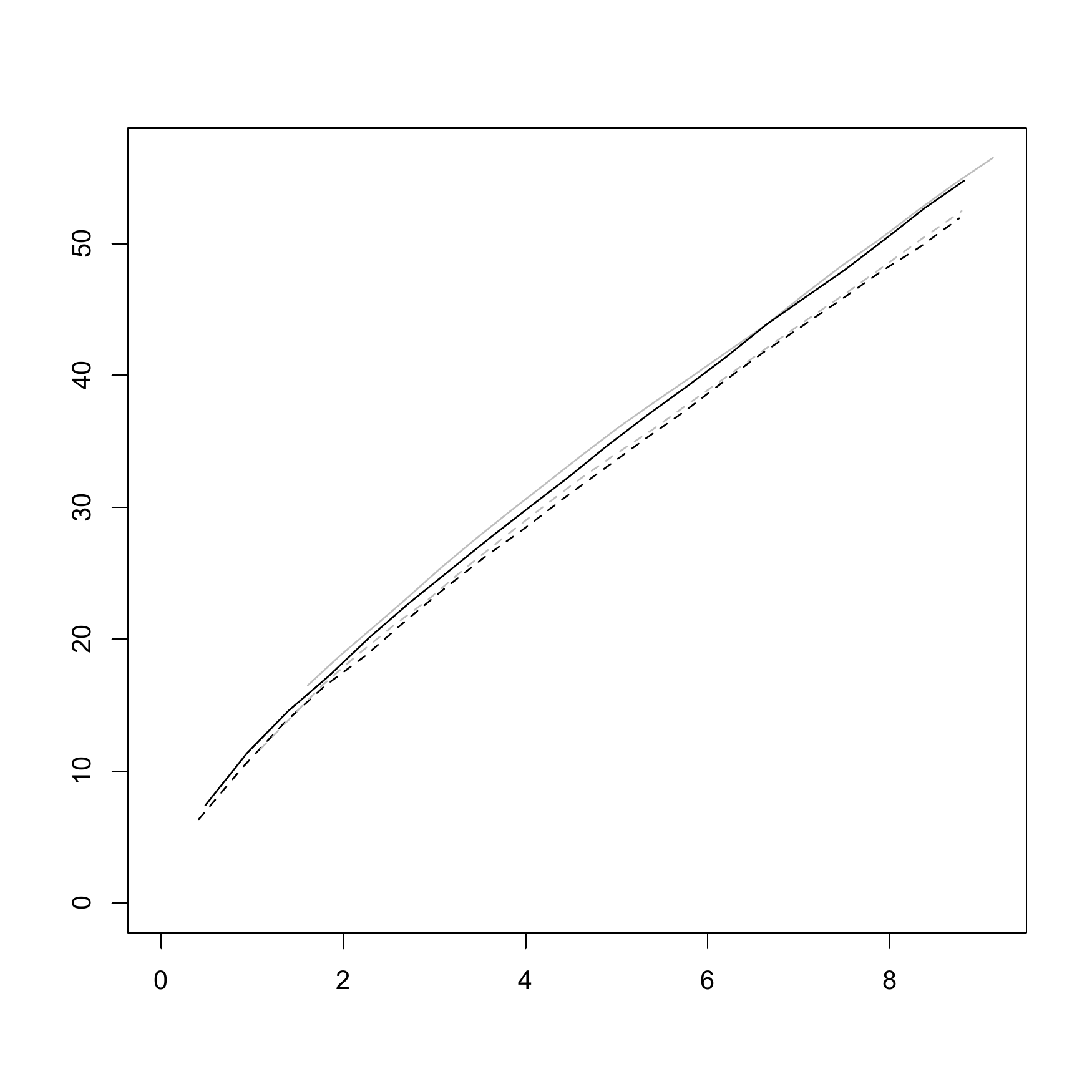}\\
$|A|=1$ & $|A|=3$ \\
\includegraphics[width=0.5\linewidth, height=0.5\linewidth]{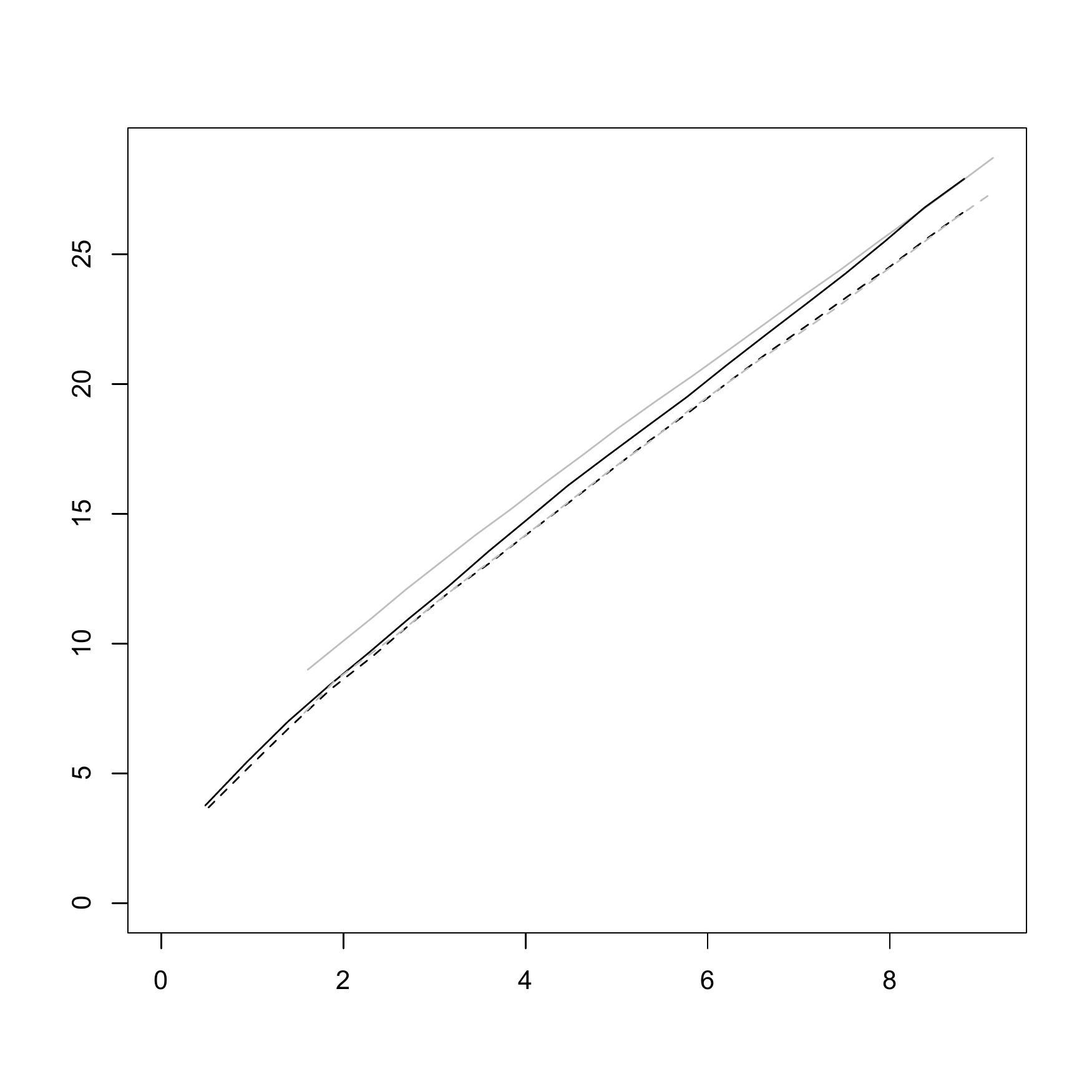}
&
\includegraphics[width=0.5\linewidth, height=0.5\linewidth]{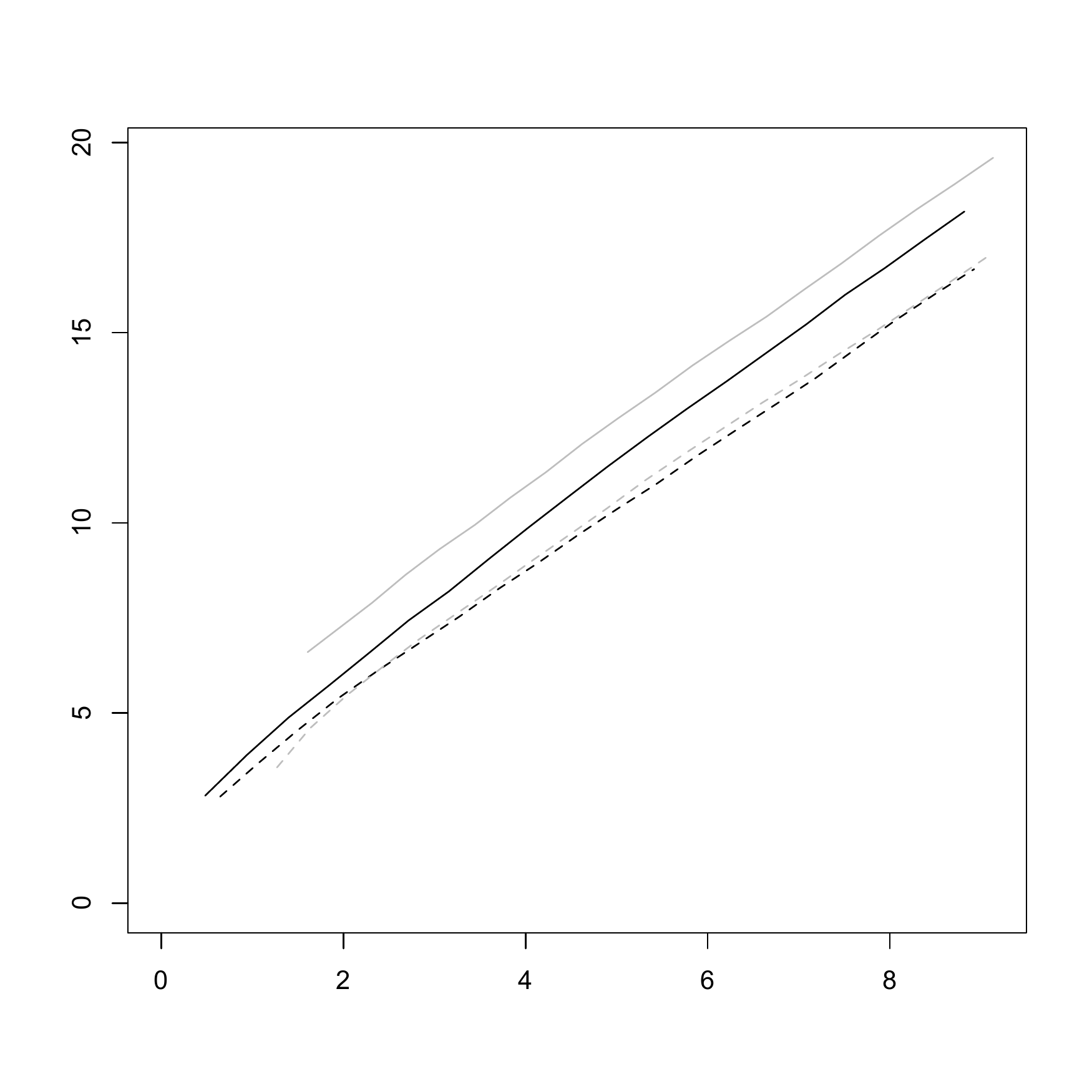}\\
$|A|=6$  & $|A|=9$ 
\end{tabular}
 \caption{ Expected sample size  against the type-I error probability in log-scale  for the   \mbox{G-SLRT} (soft lines) and the \mbox{M-SLRT} (dark lines). That is, if $\cA$ is the affected subset, we plot 
 $ \Exp^{\cA}[T]$ (vertical axis) against   $|\log \Pro_0(d=1)|$ (horizontal axis) for the following cases: $|\cA|=1,3,6,9$.  For both tests, solid lines refer to the case of no prior information ($\overline{\cP}_{K}$), whereas dashed lines refer to the case  that  the size of the affected subset is known in advance 
 ($\cP_{|\cA|}$). \label{fig1} }
\end{figure}

 For  a fair comparison between the \mbox{G-SLRT} and the \mbox{M-SLRT}, we need to compare their expected sample sizes when the two schemes have the same error probabilities.  In  Figure \ref{fig1} we plot the expected sample size of each test against the logarithm of the type-I error probability for  different cases regarding the size of the affected subset.  Specifically, if $\cA$ is the affected subset, we plot   $\Exp^{\cA}[T]$ (vertical axis) against   $|\log \Pro_0(d=1)|$ (horizontal axis) for the following cases: $|\cA|=1,3,6,9$.  The dashed lines correspond to the versions of the two schemes when the size of the affected subset is known ($\overline{\cP}_{|\cA|}$). The solid lines correspond to the versions of the two schemes with no prior information ($\overline{\cP}_K$).   The dark  lines correspond to \mbox{M-SLRT}, whereas the grey lines to \mbox{G-SLRT}.

 We observe that when we design the two tests knowing the size of the affected subset, then their performance  is essentially identical. However, when we design the two tests assuming no  prior information,  the \mbox{G-SLRT}  performs slightly better (resp. worse) than  the \mbox{M-SLRT} in the case where the signal is present in  a small (resp. large) number of channels, at least for large and moderate error probabilities. The operating characteristics of the two tests become almost identical as the type-I error goes to 0, as expected.   Note however that when  the number of affected channels is large, the signal-to-noise ratio is high. Therefore, the ``absolute" loss of the \mbox{G-SLRT}  in these cases is small.

\begin{figure}
\centering
\begin{tabular}{cc}
\includegraphics[width=0.5\linewidth, height=0.5\linewidth]{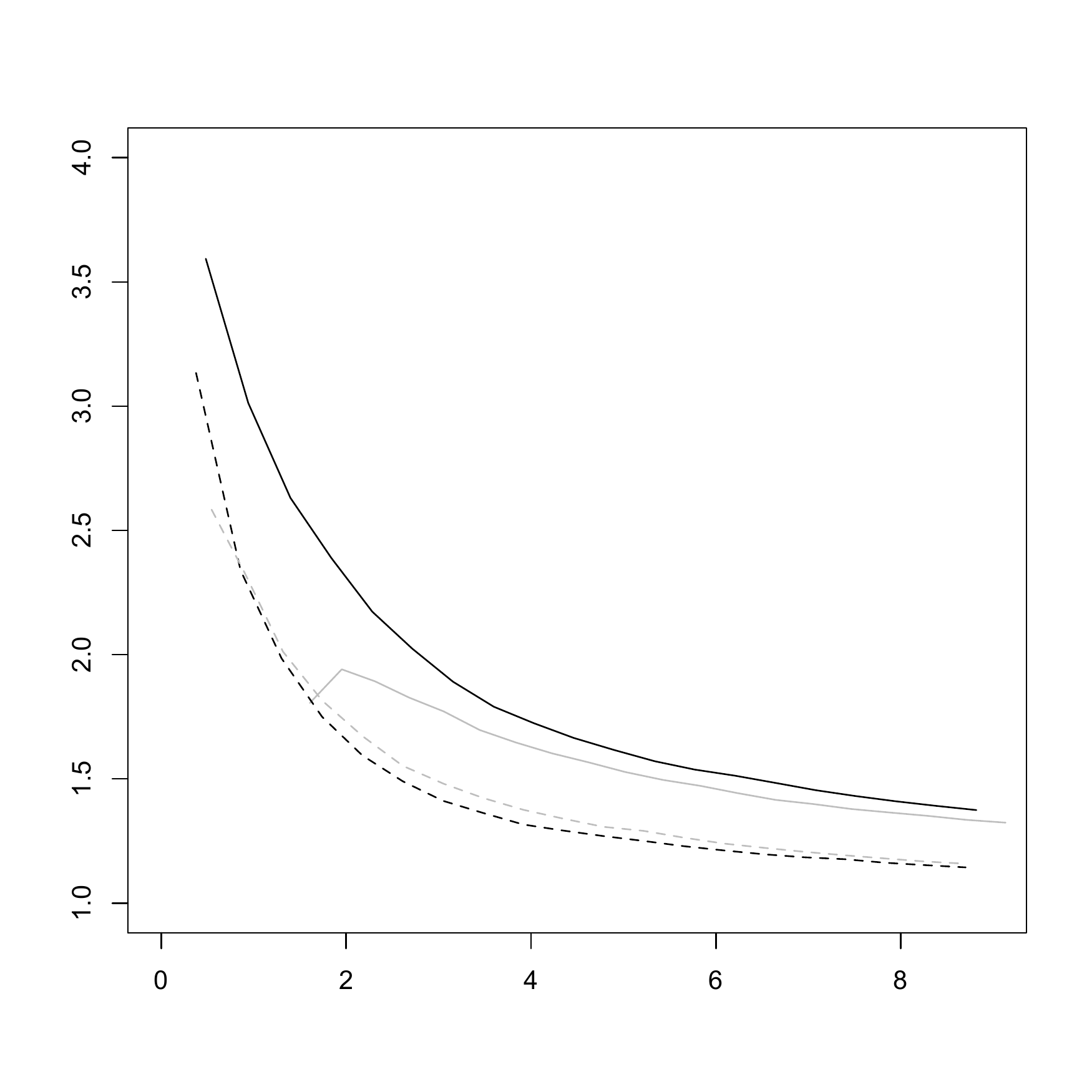}
&
\includegraphics[width=0.5\linewidth, height=0.5\linewidth]{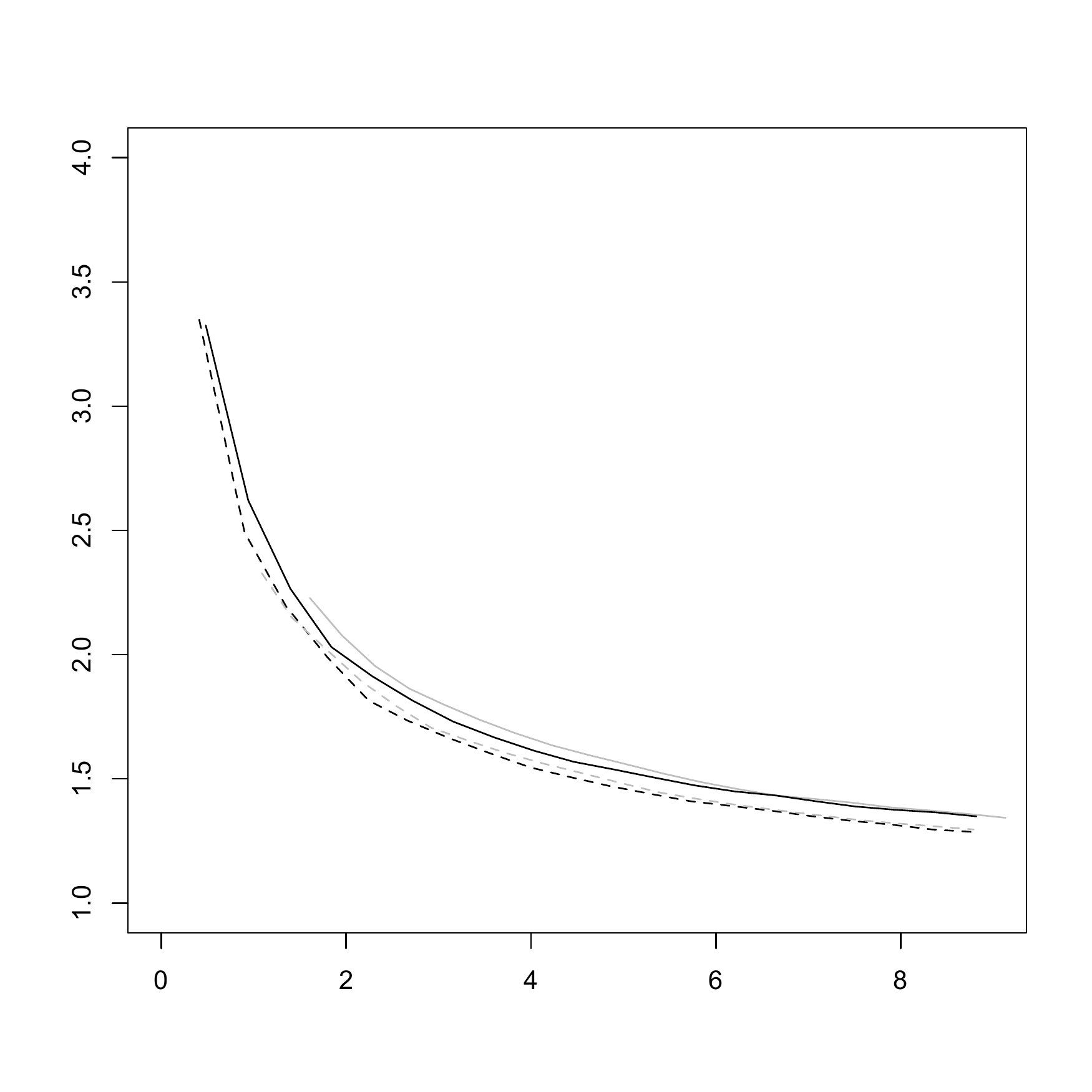}\\
$|A|=1$ & $|A|=3$ \\
\includegraphics[width=0.5\linewidth, height=0.5\linewidth]{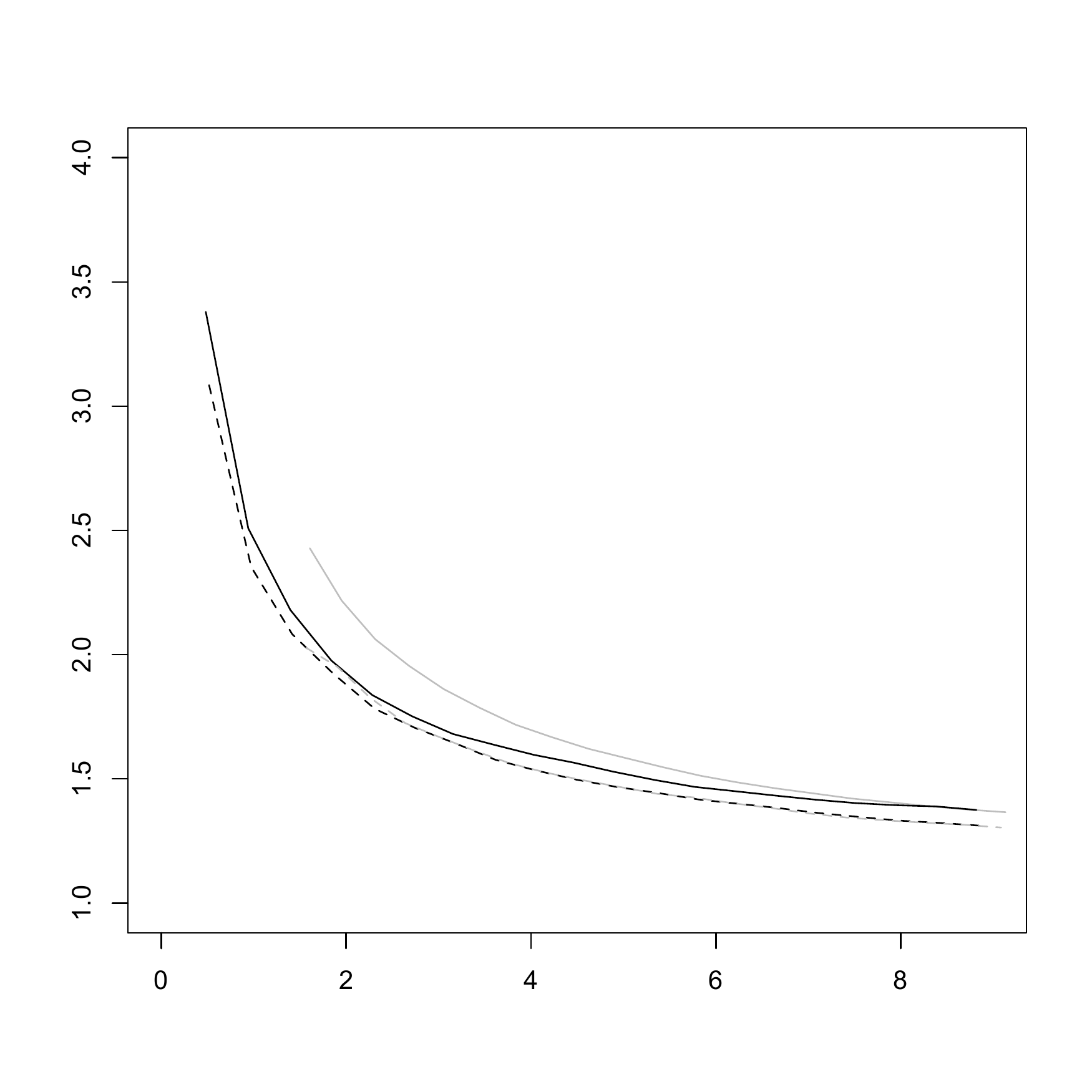}
&
\includegraphics[width=0.5\linewidth, height=0.5\linewidth]{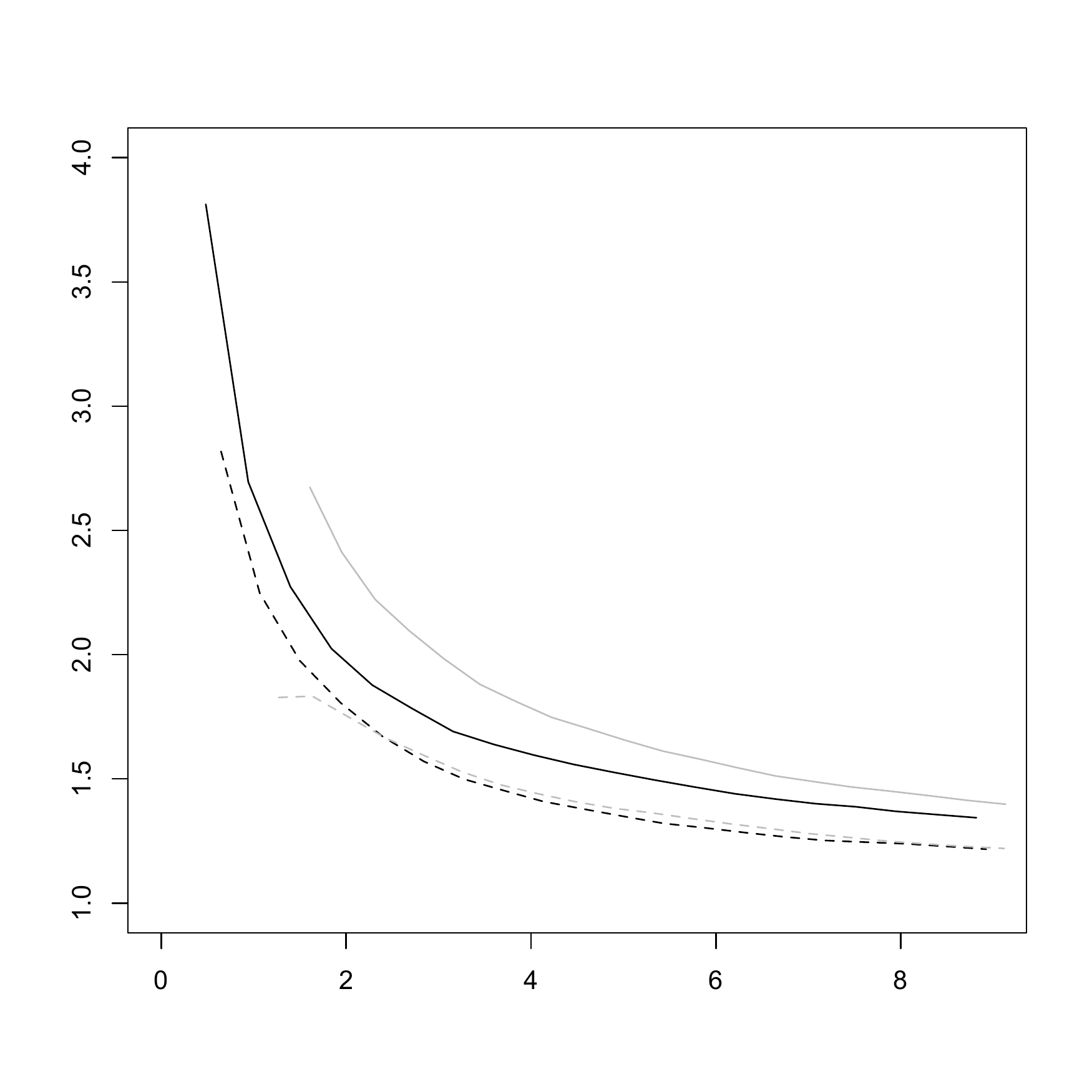}\\
$|A|=6$  & $|A|=9$ 
\end{tabular}
 \caption{ Normalized expected sample size  against the type-I error probability in log-scale  for the   \mbox{G-SLRT} (soft lines) and the \mbox{M-SLRT} (dark lines). That is, if $\cA$ is the affected subset, we plot 
 $|\cA| I_1 \Exp^{\cA}[T] /\,  | \log \Pro_0(d=1)|$ (vertical axis) against   $|\log \Pro_0(d=1)|$ (horizontal axis) for the following cases: $|\cA|=1,3,6,9$.  For both tests, solid lines refer to the case of no prior information ($\overline{\cP}_{K}$), whereas dashed lines refer to the case  that  the size of the affected subset is known in advance 
 ($\cP_{|\cA|}$). \label{fig2} }
\end{figure}

Finally, in   Figure \ref{fig2} we plot the normalized expected sample of each test under the alternative hypothesis against the logarithm of the type-I error probability for  different cases regarding the size of the affected subset.  That is, if $\cA$ is the affected subset, we plot   $|\cA| I_1 \Exp^{\cA}[T]/ |\log \Pro_0(d=1)|$ (vertical axis) against   $|\log \Pro_0(d=1)|$ (horizontal axis) for the following cases: $|\cA|=1,3,6,9$.  Again, the dashed lines correspond to the versions of the two schemes when the size of the affected subset is known ($\overline{\cP}_{|\cA|}$). The solid lines correspond to the versions of the two schemes with no prior information ($\overline{\cP}_K$).   The dark  lines correspond to the \mbox{M-SLRT}, whereas the gray lines to the \mbox{G-SLRT}.  Our asymptotic theory suggests that the curves in Figure \ref{fig2} converge to 1, and this is also  verified by our graph. The convergence is relatively slow in most cases, which can be explained by the fact that we do not normalize  the expected sample sizes by  the optimal performance, but with an asymptotic lower bound on it. 

\section{Conclusion and Remarks} \label{sec:Conclusions}
 
 We considered  the  problem of sequential detection of an unknown number of signals in multiple data streams and studied two families of sequential tests. The first,  \mbox{G-SLRT},    is based on maximizing the  likelihood ratios between the ``signal and noise'' and ``noise only'' hypotheses. The second, \mbox{M-SLRT},
 is based on a mixture (weighted sum) of likelihood ratios. Based on the concept of \textit{$r$-complete}  convergence, we developed a general theory that allows for the study of asymptotic properties of the above sequential tests for very general non-i.i.d. models without assuming any particular  structure for the observations 
 apart from an asymptotic stability property of the local log-likelihood ratios.  Specifically, under the assumption that the log-likelihood ratios in channels 
converge $r$-completely when suitably normalized,  we were able to show that both tests  asymptotically minimize moments of the sample size up to order $r$ as the probabilities of errors approach zero. Moreover, in the special case that  the local log-likelihood ratios have independent (\textit{but not necessarily identically distributed})  increments and converge only \textit{almost surely} when suitably normalized,  we showed that  both tests asymptotically minimize \textit{all} moments of the sample size. 

These asymptotic optimality results were shown under the assumption of an arbitrary class of possibly affected subsets. They are thus valid for both structured and unstructured multistream hypothesis testing problems. Moreover, we illustrated this  general sequential hypothesis testing theory using several meaningful examples including Markov and hidden Markov models, as well as a multichannel generalization of the famous invariant $t$-SPRT. Finally, when compared using a simulation study, the \mbox{G-SLRT} (\mbox{M-SLRT})  was found to perform better when a small (large) number of channels is affected and there is no prior information regarding the affected subset. On the other hand, the two procedures were found to perform similarly when the size of the affected subset is known in advance. 

When the observations in channels are i.i.d., even if they differ across channels, we can obtain  stronger and more refined results  for the proposed procedures  along the lines of our previous 
works~\cite{FellourisTartakovsky-SQA2012,FellourisTartakovsky-SS2013},   such as  near-optimality and higher order approximations. 
These results are based on nonlinear renewal theory and will be presented in the companion paper \cite{FTPart2}.  
Moreover, it is also possible to  generalize our  asymptotic analysis by allowing the number of channels to approach infinity, which is also a topic of the companion paper \cite{FTPart2}.

\section*{Acknowledgment}

We would like to thank Dr. G. V. Moustakides for stimulating discussions and for  helpful suggestions in Section~\ref{ssec:GSLRT} of the paper.


\appendix

\renewcommand{\theequation}{A.\arabic{equation}}
\setcounter{equation}{0}

  \renewcommand{\thelemma}{\thesubsection.\arabic{lemma}}

\subsection{Proofs}

\begin{lemma} \label{Lem: maxSLLN}
Let $\{Z_t\}_{t\in \bN}$ be a stochastic process defined on some probability space $(\Omega, \cF, \Pro)$ and let $\Exp$ be the corresponding expectation.  
Suppose that $t^{-1}Z_t$ converges almost surely to a finite, positive  constant $I$ as $t\to\infty$. Then
\[
\lim_{M\to\infty}\Pro\set{\frac{1}{M} \max_{1 \le t \le M} Z_t > (1+\varepsilon) I} =0 \quad \text{for all} ~ \varepsilon >0 .
\]

\end{lemma}

\begin{IEEEproof}
Write $Y_t= t^{-1} Z_t - I$, $U_N(\varepsilon)=\bigcup_{t>N}\{|Y_t| \ge \varepsilon I\}$,
\[
P_M(\varepsilon)=\Pro\set{\frac{1}{M} \max_{1 \le t \le M} Z_t > (1+\varepsilon) I} \quad \text{and} \quad  P_{M,N}(\varepsilon)= \Pro\set{\frac{1}{M}\max_{1\le t \le N} Z_t \ge (1+\varepsilon) I } .
\]
For any fixed $1\le N\le M$, by the addition rule  we have 
\begin{align*}
P_M(\varepsilon) &\le  P_{M,N}(\varepsilon) +  \Pro\set{\max_{N < t \le M} Z_t 
 \ge (1+\varepsilon) I M}   .
 \end{align*}
For the second term we have the following chain of inequalities
 \begin{align*}
 \Pro\set{\max_{N < t \le M} Z_t 
 \ge (1+\varepsilon) I M}  &\leq  \Pro\set{\max_{N < t \le M} (Z_t - I t) \ge \varepsilon I M}  \\
& =  \Pro\set{\max_{N < t \le M} t Y_t \ge \varepsilon I M}
\\
& \le   \Pro\set{\max_{N < t \le M} Y_t \ge \varepsilon I}
\le   \Pro\set{\max_{t>N} Y_t \ge \varepsilon I}
\\
& \le  \Pro\set{\max_{t>N} |Y_t| \ge \varepsilon I}  \le  \Pro( U_N(\varepsilon)) .
\end{align*}
Thus, for any $N\ge 1$, $M\geq N$  and $\varepsilon >0$ we have 
\begin{align} \label{NM}
P_M(\varepsilon) &\le  P_{M,N}(\varepsilon) +  \Pro(U_N(\varepsilon)).
 \end{align}
 Since $\Pro(|Z_t|<\infty)=1 $ for every $t \in \bN$, from Markov's inequality it follows  that for any $N\ge 1$  and $\varepsilon >0$ we have 
$
\lim_{M\to\infty} P_{M,N}(\varepsilon)=0
$
 and, consequently, letting $M \rightarrow \infty$ in \eqref{NM} we obtain
\begin{equation}\label{limsupP1M}
\limsup_{M\to\infty} P_{M}(\varepsilon) \le  \Pro (U_N (\varepsilon)) .
\end{equation}
But from  the definition of a.s.\ convergence and the assumption of the lemma it follows that,
for any $\varepsilon>0$, $\lim_{N\to\infty} \Pro \set{U_{N}(\varepsilon)} =0$.
Hence, letting  $N\to\infty$ in \eqref{limsupP1M}, we obtain the assertion of the lemma. 
\end{IEEEproof}

\begin{IEEEproof}[Proof of Lemma \ref{rcompleteglobal}]
Consider an arbitrary subset $\cA \in \cP$ and $\varepsilon>0$.  We have to show that 
\begin{equation} \label{rcompleteglobaldetailed}
\begin{aligned}
 \sum_{t=1}^\infty t^{r-1} \Pro^\cA\brc{\abs{\frac{1}{t} Z_{t}^{\cA} -I_1^\cA} > \varepsilon} < \infty , \quad  \sum_{t=1}^\infty t^{r-1} \Pro_0 \brc{\abs{\frac{1}{t} Z_{t}^{\cA} +I_0^\cA} > \varepsilon} < \infty
 \end{aligned}
 \end{equation}
 whenever $r$-complete conditions \eqref{rcompleteLLRsdetailed} for $t^{-1} Z_t^k$ hold for all $k=1,\dots,K$.

For every  $t \in \bN$, we have $|t^{-1}Z_t^\cA - I_1^\cA| \leq  \sum_{k \in \cA}  |t^{-1}Z_t^k - I_1^k|$,  and therefore,
\begin{equation} \label{useful}
\set{t\ge 1: |t^{-1} Z_t^\cA - I_1^\cA| > \varepsilon} \subset \bigcup_{k \in \cA} \set{t\ge 1: |t^{-1}Z_t^k - I_1^k| > \varepsilon/|\cA|} .
\end{equation}
Hence,
$$
\Pro^\cA\brc{\abs{\frac{1}{t} Z_{t}^{\cA} -I_1^\cA} > \varepsilon} \leq \sum_{k \in \cA} \Pro_{1}^{k} \brc{\abs{ \frac{1}{t}Z_t^k - I_1^k} > \frac{\varepsilon}{|\cA|} }
 $$
 and, consequently,   by  \eqref{rcompleteLLRsdetailed},
 $$
 \sum_{t=1}^\infty t^{r-1} \Pro^\cA\brc{\abs{\frac{1}{t} Z_{t}^{\cA} -I_1^\cA} > \varepsilon} \leq   \sum_{k \in \cA} \sum_{t=1}^\infty t^{r-1}  \Pro_{1}^{k} \brc{\abs{ \frac{1}{t} Z_t^k - I_1^k} > \frac{\varepsilon}{|\cA|} }<\infty.
 $$
 The proof of the first inequality in \eqref{rcompleteglobaldetailed} is essentially similar. \\
\end{IEEEproof}

\begin{IEEEproof}[Proof of Lemma \ref{independent}]
Let $\nu_k(b) := \inf\{t: S_t^k >b\}$.  It is clear that 
$\nu(b) \geq  \nu_k(b)$.  From the SLLN it follows that  $\nu_k (b)$ is almost surely finite for any given $b>0$  and   $\nu_k (b) \rightarrow \infty$ almost surely  as $b \rightarrow \infty$. Then,  with probability 1 we have  $S_{\nu_k (b)} \geq b$ and 
$$
\frac{\nu(b)}{b}  \geq \frac{\nu_k (b)}{b} \geq \frac{\nu_k (b)}{S_{\nu_k (b)}}  \underset{b \rightarrow \infty} \longrightarrow \frac{1}{\mu_k}.
$$
Since this is true for any  $k$, we obtain
$$
\liminf_{b \rightarrow \infty} \frac{\nu(b)}{b}  \geq \left( \min_{1 \leq k \leq K} \mu_k \right)^{-1}.
$$
In order to prove the reverse inequality, we observe that
$$
\sum_{k=1}^{K}  S^k_{\nu (b)} \ind{ \nu(b)=\nu_{k}(b)} \leq b + \sum_{k=1}^{K}  \xi^k_{\nu (b)} \ind{ \nu(b)=\nu_{k}(b)}  ,
$$
since   for every $k$ we have $S^k_{\nu_k (b)} \leq b+ \xi^k_{\nu_k (b)}.$
Consequently, 
$$
\min_{1 \leq k \leq K}   S^k_{\nu (b)}  \leq b +  \max_{1 \leq k \leq K} \xi^k_{\nu(b)}. 
 $$
and 
 $$
\min_{1 \leq k \leq K}   \frac{S^k_{\nu (b)}}{\nu(b)}  \leq \frac{b}{\nu(b)} +  \max_{1 \leq k \leq K}  \frac{\xi^k_{\nu(b)}}{\nu(b)},
 $$
 which implies that 
  $$
\liminf_{b \rightarrow \infty}  \frac{b}{\nu(b)} \geq  \min_{1 \leq k \leq K}  \mu_k,
$$
since 
$$
\frac{\xi^k_t}{t}=\frac{S^k_t}{t}- \frac{t-1}{t} \frac{S^k_{t-1}}{t-1} \rightarrow 0.
$$
It remains to show that $(\nu_b/b)_{b>0}^r$ is uniformly integrable for every $r>0$ when \eqref{condition} holds.  It suffices to restrict ourselves to $b \in \bN$.  Similarly to \cite[Theorem 2.5.1, p. 57]{Gut-book88},  we observe that for any $b,c \in \bN$ we have  $\nu(b+c) \leq \nu(b)+ \nu(c;b),$ where 
$$
\nu(c;b):=\inf\left\{t> \nu(b):  S_t^k -S^k_{\nu(b)} >c \quad  \forall \; 1\leq k \leq K \right\}.
$$
By induction,
$$
\nu(b) \leq \sum_{n=0}^{b-1} \nu(1;n)
$$
and, consequently,
$$
||\nu(b)||_r \leq \sum_{n=0}^{b-1} ||\nu(1;n)||_r \leq b \sup_{n \in \bN}  ||\nu(1;n)||_r 
$$
and 
$$
||\nu(b) / b ||_r \leq  \sup_{n \in \bN}  ||\nu(1;n)||_r .
$$
It remains to show that the upper bound is finite when \eqref{condition} holds.  Indeed,  for any $m \in \bN$,
\begin{align*}
  \Pro( \nu(1;n) >m) &=\Pro \left( \max_{ \nu(n) <t \leq  \nu(n)+m}  (S_t^k -S^k_{\nu(n)}) \leq 1 \quad \text{for some $1 \leq  k \leq K$} \right)  \\
  &\leq  \sum_{k=1}^{K} \Pro\left(  \max_{\nu(n) <t \leq  \nu(n)+m}  (S_t^k -S^k_{\nu(n)}) \leq 1 \right) \\
    &\leq  \sum_{k=1}^{K}   \Pro \left(  S^k_{\nu(n)+m} -S^k_{\nu(n)} \leq 1\right)
\end{align*} 
and  from Markov's inequality we obtain for any $\lambda \in (0,1)$:
\begin{align*}
 \Pro(  S^k_{\nu(n)+m} -S^k_{\nu(n)} \leq 1)  
    &\leq \Pro\left( \exp \left\{ -\lambda ( S^k_{\nu(n)+m} -S^k_{\nu(n)} \right\} \geq e^{-\lambda} \right)  \\
  &\leq   e^{\lambda}  \, \Exp\left[ \exp\left\{ -\lambda ( S^k_{\nu(n)+m} -S^k_{\nu(n)} ) \right\} \right] \\
    &\leq   e^{\lambda} \,  \Exp\left[  \prod_{u=\nu(n)+1}^{\nu(n)+m}\exp\left\{ -\lambda \xi^k_u \right\} \right].
   \end{align*} 
If we set $\beta_{k}(\lambda):= \sup_{n \in \bN} \Exp\left[ \exp\{\lambda (\xi^k_n)^{-} \} \right]$, then  from Lemma~\ref{prod_lem} (see below) we have 
\begin{align*}
 \Exp\left[  \prod_{u=\nu(n)+1}^{\nu(n)+m}\exp\left\{ -\lambda \xi^k_u \right\} \Big| \nu(n) \right] 
          &=   \prod_{u=\nu(n)+1}^{\nu(n)+m}  \Exp\left[ \exp\{-\lambda \xi^k_u\} \right]   \\
        &\leq   \prod_{u=\nu(n)+1}^{\nu(n)+m}  \Exp\left[ \exp\{\lambda (\xi^k_u)^{-} \}  \right]   \leq   \beta_{k}(\lambda) ^m.
   \end{align*}
We conclude that 
$$
   \Pro( \nu(1;n) >m) \leq  e^{\lambda} \sum_{k=1}^{K}  \beta_{k}(\lambda)  ^m \leq (K e^{\lambda})  \left(\max_{1 \leq k \leq K} \beta_{k}(\lambda) \right)^m,
$$
which implies that $\sup_{n \in \bN}  ||\nu(1;n)||_r <\infty$ for every $r >0$ and  completes the proof. 
\end{IEEEproof}

\begin{lemma} \label{prod_lem}
Let $\xi=(\xi_t)_{t \in \bN}$ be a  sequence of  positive, independent random variables  on some probability space $(\Omega, \cF, \cP)$. Suppose that  $\Exp[\xi_t]<\infty$ for every $t \in \bN$, where 
$\Exp$ is expectation with respect to $\Pro$. Let $T$ be a stopping time with respect to the filtration generated by $\xi$. Then, for every deterministic integer $m \in \bN$ we have 
\begin{equation} \label{prod}
 \Exp \left[ \prod_{u=T+1}^{T+m} \xi_u  \, \Big|  \, T\right] =  \prod_{u=T+1}^{T+m} \Exp[\xi_u].
\end{equation}
When in particular $\Exp[\xi_t] \leq c$ for every $t \in \bN$ for some constant $c$, then 
$$
 \Exp \left[ \prod_{u=T+1}^{T+m} \xi_u  \right] \leq c^m.
$$
\end{lemma}
\begin{IEEEproof}
For any $t \in \bN$ we have 
\begin{align*}
\Pro(T=t)  \, \Exp \left[ \prod_{u=T+1}^{T+m} \xi_u  \, \Big|  \, T=t\right]  &= \Exp \left[ \prod_{u=t+1}^{t+m} \xi_u \;  ; \;  T=t\right] 
\\ &= \Pro(T=t) \; \Exp \left[ \prod_{u=t+1}^{t+m} \xi_u  \right]  
=  \Pro(T=t)  \prod_{u=t+1}^{t+m}  \Exp[ \xi_u],
\end{align*}
where the second equality holds because the random variables  $\{\xi_u, t+1\leq u \leq t+m\}$  are independent of the event $\{T=t\}$, which  depends on  $\{\xi_u, 1\leq u \leq t\}$. This proves  \eqref{prod}. 
 \end{IEEEproof}

\subsection{Details on the AR model} \label{appen;AR}

Here, we  provide  more details regarding the proof of the $r$-complete convergence in the autoregressive model of Subsection \ref{ssec:AR}.  We essentially need to show that conditions $(C_1)$ and $(C_2)$ in \cite[Sec 5]{PergTar} hold.  Define 
\[
\hat{g}_k(x) := \int_{-\infty}^\infty g_k(y,x) \varphi(y-\rho_1^k x) dy = \frac{(\rho_1^k-\rho_0^k)^2 x^2}{2} .
\]
We have
\begin{equation}\label{sec:Ex.2-1}
\sup_{y,x\in(-\infty,\infty)} \frac{\vert g_k(y,x)\vert}{1+|y|^2+|x|^2} \le Q
\quad\mbox{and}\quad \sup_{x\in(-\infty,\infty)} \frac{\hat{g}(x)}{1+|x|^2} \le Q,
\end{equation}
where 
$$
Q=\max\set{1,\frac{\vert(\rho_1^k)^{2}-(\rho_0^k)^{2}\vert+(\rho_1^k-\rho_0^k)^{2}+1}{2}}.
$$ 
Define also the Lyapunov function $V(x)= Q(1+|x|^2)$. Obviously, 
$$
\lim_{|x|\to\infty} \frac{\Exp_{x,1}^k[ V(X_1^k)]}{V(x)} =\lim_{|x|\to\infty} \frac{1+\Exp[ |\rho_1^k x+\xi_{1}^k|^2]}{1+|x|^2} =|\rho_1^k|^2<1,
$$
where $\Exp_{x,1}^k$ stands for expectation under $\Pro_{x,1}^k=\Pro_1^k(\cdot | X_0^k=x)$. Therefore, for any $|\rho_1^k|^2<\varrho<1$ there exist $D>0$ such that the condition $(C_1)$ in \cite[Sec 5]{PergTar} holds 
with $C=[-n, n]$ for every $n\ge1$. 

Next, since all the moments of $\xi_1^k$ are finite, it follows that $\Exp[\vert w_{0}^k] \vert^{r}<\infty$ and 
$\Exp [\vert w_1^k] \vert^{r}<\infty$ for all $r \ge 1$. Moreover, taking into account the ergodicity properties, we obtain that for any $x\in(-\infty,\infty)$
\begin{equation}\label{sec:Ex.2-3n}
\lim_{t\to\infty} \Exp_{x,0}^k [\vert X_t^k\vert^{r}]=\Exp[\vert w_0^k\vert^{r}]<\infty \quad\mbox{and}\quad
\lim_{t\to\infty}\Exp_{x,1}^k [\vert X_t^k\vert^{r}]=\Exp [\vert w_1^k \vert^{r}] <\infty.
\end{equation}
Observe that under $\Pro_{x,1}^k$ for any $t\ge 1$
$$
X_t^k=(\rho_1^k)^t x +\sum^{t}_{\ell=1}( \rho_1^k)^{t-\ell} \xi_\ell^k.
$$
Hence, for any $r \ge 1$,
$$
\Exp_{x,1}^k [\vert X_t^k \vert^{r}] \le 2^{r} \left(\vert x\vert^{r} + \Exp_{0,1}^k \vert X_t^k \vert^{r}\right),
$$
i.e., using the last convergence in \eqref{sec:Ex.2-3n} we obtain that for some $C^{*}>0$
$$
M^{*}(x)=\sup_{t\ge 1} \Exp_{x,1}^k [\vert X_t^k \vert^{r} ] \le C^{*} (1+\vert x\vert^{r}).
$$
Using now the first convergence in \eqref{sec:Ex.2-3n} we obtain that $\sup_{t\ge 1}\Exp_1^k [M^{*}(X_t^k)]<\infty$.
So, the upper bounds in \eqref{sec:Ex.2-1} imply the condition ($C_2$) in \cite[Sec 5]{PergTar}. 

\bibliographystyle{IEEEtran}
%
\bibliography{kniga_plus_IEEEIT}


\end{document}